\documentclass[twocolumn,amssymb,aps,nofootinbib]{revtex4}
\usepackage{times,mathptmx,amsthm,amsmath}
\usepackage{pstricks,pst-node,graphics}
\makeatletter

\renewcommand{\p@subsection}{}

\renewcommand{\p@subsubsection}{}

\makeatother

\newtheorem{theorem}{Theorem}
\newtheorem{lemma}[theorem]{Lemma}

\newtheorem{question}[theorem]{Question}
\theoremstyle{remark}
\newtheorem*{remark}{Remark}

\newcommand{\V}{\mathrm{V}}     \newcommand{\VH}{\mathrm{VH}}
\newcommand{\VHP}{\mathrm{VHP}} \newcommand{\HT}{\mathrm{HT}}
\newcommand{\QT}{\mathrm{QT}}   \newcommand{\U}{\mathrm{U}}
\newcommand{\UU}{\mathrm{UU}}   \renewcommand{\O}{\mathrm{O}}
\newcommand{\OO}{\mathrm{OO}}   \newcommand{\UO}{\mathrm{UO}}

\newcommand{\ba}{\bar{a}} \newcommand{\bb}{\bar{b}}
\newcommand{\bc}{\bar{c}} \newcommand{\bs}{\bar{s}}
\newcommand{\bt}{\bar{t}} \newcommand{\bx}{\bar{x}}
\newcommand{\by}{\bar{y}} \newcommand{\bz}{\bar{z}}
\newcommand{\ii}{\mathbf{i}}
\newcommand{\tA}{\widetilde{A}}
\renewcommand{\sl}{\mathfrak{sl}}
\DeclareMathOperator{\Pf}{Pf}

\newcommand{\vx}{\vec{x}} \newcommand{\vy}{\vec{y}}
\newcommand{\vq}{\vec{q}} \newcommand{\vone}{\vec{1}}

\newcommand{\eatline}{\vspace{-\baselineskip}}

\newcommand{\Ie}{\textit{I.e.}}

\newenvironment{fullfigure}[2]
    {\begin{figure}[htb]\begin{center}\def\ffa{#1}\def\ffb{#2}}
    {\vspace{\baselineskip}\caption{\ffb.}\label{\ffa}\end{center}\end{figure}}

\newenvironment{fulltable}[2]
    {\begin{table}[htb]\begin{center}\def\ffa{#1}\def\ffb{#2}}
    {\vspace{\baselineskip}\caption{\ffb.}\label{\ffa}\end{center}\end{table}}

\newcommand{\tab}[1]{Table~\ref{#1}}
\newcommand{\fig}[1]{Figure~\ref{#1}}
\newcommand{\thm}[1]{Theorem~\ref{#1}}

\psset{linewidth=.5pt,dash=4pt 4pt,unit=.25in,arrowsize=2pt 3}
\newgray{gray40}{.6}
\SpecialCoor

\newcommand{\mto}{\lput{:U}{\pspicture(0,0)(0,0)
\psline[arrows=->](2.3pt,0)(2.4pt,0)\endpspicture}}
\newcommand{\mfro}{\lput{:U}{\pspicture(0,0)(0,0)
\psline[arrows=->](-2.3pt,0)(-2.4pt,0)\endpspicture}}

\newcommand{\verta}[1]{\pspicture[.6](-.7,-1.9)(.7,.7)
\pcline(0,0)(0,.7)\mto    \pcline(0,0)(.7,0)\mfro
\pcline(0,0)(-.7,0)\mfro  \pcline(0,0)(0,-.7)\mto
\rput(0,-1.2){$#1$}\endpspicture}
\newcommand{\vertb}[1]{\pspicture[.6](-.7,-1.9)(.7,.7)
\pcline(0,0)(0,.7)\mfro   \pcline(0,0)(.7,0)\mto
\pcline(0,0)(-.7,0)\mto   \pcline(0,0)(0,-.7)\mfro
\rput(0,-1.2){$#1$}\endpspicture}
\newcommand{\vertc}[1]{\pspicture[.6](-.7,-1.9)(.7,.7)
\pcline(0,0)(0,.7)\mto    \pcline(0,0)(.7,0)\mto
\pcline(0,0)(-.7,0)\mfro  \pcline(0,0)(0,-.7)\mfro
\rput(0,-1.2){$#1$}\endpspicture}
\newcommand{\vertd}[1]{\pspicture[.6](-.7,-1.9)(.7,.7)
\pcline(0,0)(0,.7)\mfro   \pcline(0,0)(.7,0)\mfro
\pcline(0,0)(-.7,0)\mto   \pcline(0,0)(0,-.7)\mto
\rput(0,-1.2){$#1$}\endpspicture}
\newcommand{\verte}[1]{\pspicture[.6](-.7,-1.9)(.7,.7)
\pcline(0,0)(0,.7)\mto    \pcline(0,0)(.7,0)\mfro
\pcline(0,0)(-.7,0)\mto   \pcline(0,0)(0,-.7)\mfro
\rput(0,-1.2){$#1$}\endpspicture}
\newcommand{\vertf}[1]{\pspicture[.6](-.7,-1.9)(.7,.7)
\pcline(0,0)(0,.7)\mfro   \pcline(0,0)(.7,0)\mto
\pcline(0,0)(-.7,0)\mfro  \pcline(0,0)(0,-.7)\mto
\rput(0,-1.2){$#1$}\endpspicture}

\begin{document}
\title{Symmetry classes of alternating-sign matrices under one roof}
\author{Greg Kuperberg}
\email{greg@math.ucdavis.edu}
\thanks{Supported by NSF grants DMS \#9704125 and DMS \#0072342, and by a
    Sloan Foundation Research Fellowship.}
\affiliation{UC Davis}
\begin{abstract}
In a previous article \cite{Kuperberg:asm}, we derived the alternating-sign
matrix (ASM) theorem from the Izergin-Korepin determinant
\cite{Izergin:finite,KBI:qism,ICK:determinant} for a partition function for
square ice with domain wall boundary.  Here we show that the same argument
enumerates three other symmetry classes of alternating-sign matrices:  VSASMs
(vertically symmetric ASMs), even HTSASMs (half-turn-symmetric ASMs), and even
QTSASMs (quarter-turn-symmetric ASMs). The VSASM enumeration was conjectured by
Mills; the others by Robbins \cite{Robbins:symmetry}.  We introduce several new
types of ASMs:   UASMs (ASMs with a U-turn side), UUASMs (two U-turn sides),
OSASMs (off-diagonally symmetric ASMs), OOSASMs (off-diagonally,
off-antidiagonally symmetric), and UOSASMs (off-diagonally symmetric with
U-turn sides). UASMs generalize VSASMs, while UUASMs generalize VHSASMs
(vertically and horizontally symmetric ASMs) and another new class, VHPASMs
(vertically and horizontally perverse). OSASMs, OOSASMs, and UOSASMs are
related to the remaining symmetry classes of ASMs, namely DSASMs (diagonally
symmetric), DASASMs (diagonally, anti-diagonally symmetric), and TSASMs
(totally symmetric ASMs).  We enumerate several of these new classes, and we
provide several 2-enumerations and 3-enumerations.

Our main technical tool is a set of multi-parameter determinant and Pfaffian
formulas generalizing the Izergin-Korepin determinant for ASMs and the Tsuchiya
determinant for UASMs \cite{Tsuchiya:determinant}.  We evaluate specializations
of the determinants and Pfaffians using the factor exhaustion method.
\end{abstract}
\maketitle

\section{Introduction}
\label{s:intro}

An \emph{alternating-sign matrix} (or ASM) is a matrix with entries 1,
0, and $-1$, such that the non-zero entries alternate in sign in each row and
column, and such that the first and last non-zero entry in each row and column
is 1.  Mills, Robbins, and Rumsey \cite{MRR:asm} conjectured a formula for the
number of ASMs of order $n$.  This formula was first proved by Zeilberger
in 1995 \cite{Zeilberger:asm}, and later the author found a different proof
\cite{Kuperberg:asm}:

\begin{theorem}[Zeilberger] There are
$$A(n) = \frac{1!4!7!\cdots(3n-2)!}{n!(n+1)!(n+2)!\cdots(2n-1)!}$$
$n \times n$ ASMs.
\label{th:asm}
\end{theorem}

\thm{th:asm} is part of a larger, unfinished structure in enumerative
combinatorics, much of it conjectured by Robbins
\cite{Robbins:symmetry,Robbins:story}.  The structure includes two types of
relations between alternating-sign matrices and another class of combinatorial
objects, plane partitions in boxes.  (A plane partition in a box is an order
ideal in the poset $[1,\ldots,a]\times[1,\ldots,b]\times[1,\ldots,c]$.  They
can be interpreted as a basis for the irreducible Weyl representation
$V(c\lambda_a)$ of the Lie algebra $\sl(a+b)$.)  One relation is by analogy:
The number of plane partitions in a given box is \emph{round} (meaning a
product of small factors, also called \emph{smooth}), and so is the number in
any given symmetry class.  It is usually easy to conjecture an explicit product
formula for round numbers.  Likewise Robbins found that the number of ASMs in
most (but not all!) symmetry classes also seems to be round.  The other relation
is equinumeration.  Robbins also found that there are the same numbers, both
round and not, of many types of ASMs as there are other types of plane
partitions.  To begin with there are exactly as many ASMs with no symmetry as
there are plane partitions with full symmetry, and this is what Zeilberger
more directly proved.

Our proof of \thm{th:asm} sheds no light on plane partitions, but it does rely
on a connection to another important structure in quantum algebra and
statistical mechanics, the Yang-Baxter equation.  Using the Yang-Baxter
equation, Izergin and Korepin found a determinant formula for the partition
function of square ice  with domain wall boundary conditions 
\cite{Izergin:finite,ICK:determinant,KBI:qism}.  We noted that this state model
is equivalent to certain weighted enumerations of ASMs.  Although the
determinant is singular at the point where all weights are equal, it is
generically non-singular and round along a special curve of weights with a
coordinate $q$ such that $q=1$ is the equal-weight point. (As defined in
Section~\ref{s:factor} roundness of a polynomial in $q$ is stronger than
smoothness.)  Finally the $q$-specialization of the Izergin-Korepin determinant
independently generalizes to a round determinant in two parameters $p$ and $q$.
The two-parameter determinant can be evaluated by factor exhaustion in $p$.

In this article we generalize this argument to some of the previously known
classes of ASMs and also some new ones.  The Izergin-Korepin determinant may
have seemed accidental, but we find a similar formula for each symmetry class
of ASMs after modifying the boundary conditions at the symmetry lines
(\thm{th:z}). (One case was found previously by Tsuchiya
\cite{Tsuchiya:determinant}.)  The round $q$-specialization and its
two-parameter generalization may also have seemed accidental, but we find many
such specializations complemented by four two-parameter determinants and two
three-parameter Pfaffians (Sections \ref{s:factor} and \ref{s:prod}). Besides
ordinary enumeration, the $q$-specializations also include several
$x$-enumerations (in which each orbit of $-1$ entries in the ASM has weight
$x$) with $x=2$ and $x=3$.  (Along the way we will correct an error in the
$3$-enumeration of ASMs in Reference \citealp{Kuperberg:asm} found by Robin
Chapman.)  Besides the $q$-specializations we establish that  many of the
$x$-enumerations divide each other or otherwise share large factors
(Section~\ref{s:prod}).

We speculate that our constructions are part of a yet larger and possibly less
exotic structure in quantum algebra.  In particular the solution to the
Yang-Baxter equation that we use corresponds to the $2$-dimensional
representation of the Lie algebra $\sl(2)$.  We have not investigated what
happens when $\sl(2)$ is replaced by another Lie algebra or the $2$-dimensional
representation by another representation.

\acknowledgments

The present work began with the mistake found in
Reference~\citealp{Kuperberg:asm} by Robin Chapman and with the Tsuchiya
determinant, which is the UASM case of \thm{th:z} and which was brought to the
author's attention by Jim Propp.  We would like to thank Vladimir Korepin,
Robin Chapman, and Jim Propp more generally for their attention to the author's
work.  We would also like to thank Jan de Gier, Christian Krattenthaler, Soichi
Okada, Neil Sloane, and Paul Zinn-Justin for their interest and for finding
mistakes in earlier drafts. Finally we would like to acknowledge works by
Bressoud \cite{Bressoud:proofs}, Bressoud and Propp \cite{BP:solved}, Robbins
\cite{Robbins:story}, and Zeilberger \cite{Zeilberger:asm} for spurring the
author's interest in alternating-sign matrices.

Mathematical experiments in Maple \cite{maple} were essential at every stage
of this work.  This article is typeset using REVTeX 4 \cite{revtex} and
PSTricks \cite{pstricks}.

\section{Statement of results}
\label{s:results}

In general the \emph{$x$-enumeration} of a class of ASMs is defined as the
total $x$-weight.  The $x$-weight of an ASM is $x^n$ if it has $n$ symmetry
orbits of negative entries.  \fig{f:vsasm} shows an example.  (In the figures,
we use $+$ and $-$ for $1$ and $-1$.) Note that the $2$-enumerations of any
class of ASMs is a special case (called the free fermion point in statistical
mechanics) that can often be established by other methods
\cite{Kuperberg:eklp1,Kuperberg:eklp2,Robbins:symmetry,Ciucu:cellular,%
Ciucu:reflective,JP:quartered,Jockusch:perfect,RR:asm,MRR:asm,Propp:shuffling,%
Kuo:new,Kuperberg:coker,SZ:dimer}. The reader can also consider the elementary
case of $0$-enumeration since it is also often round.

\begin{fullfigure}{f:vsasm}{A VSASM with $x$-weight $x^2$}
$$\begin{pmatrix} 0&0&0&+&0&0&0 \\ 0&0&+&-&+&0&0 \\ +&0&-&+&-&0&+ \\
0&0&+&-&+&0&0 \\ 0&+&-&+&-&+&0 \\ 0&0&+&-&+&0&0 \\ 0&0&0&+&0&0&0
\end{pmatrix}$$ \eatline
\end{fullfigure}

We summarize the conventional symmetry classes of ASMs and what is known and
conjectured about their basic enumerations.

\begin{description}
\item[$\bullet$] ASMs - alternating-sign matrices.  The $1$-, $2$-, and 
$3$-enumerations are all previously known.

\item[$\bullet$] VSASMs - vertically symmetric ASMs.  The $2$-enumeration is
previously known.  Mills conjectured the $1$-enumeration \cite{Robbins:story}.
We establish the $1$-, $2$-, and $3$-enumerations.

\item[$\bullet$] HTSASMs - half-turn symmetric ASMs.  Robbins conjectured the
$1$-enumeration and established the $2$-enumeration \cite{Robbins:symmetry}.
We establish the $1$-, $2$-, and $3$-enumerations, but only for even order.

\item[$\bullet$] QTASMs - quarter-turn symmetric ASMs.  Robbins conjectured the
enumeration.  We establish the $1$- and $2$-enumerations only for even order.
We also prove a formula for a factor of the $3$-enumeration; the other factor
does not appear to be round.

\item[$\bullet$] VHSASMs - vertically and horizontally symmetric ASMs.  Robbins
conjectured the $1$-enumeration and the $2$-enumeration is previously known
\cite{JP:quartered}.  We only establish a determinant formula.

\item[$\bullet$] DSASMs - diagonally symmetric ASMs.  Their number does not appear
to be round and no determinant or Pfaffian formula is known.

\item[$\bullet$] DASASMs - diagonally and antidiagonally symmetric ASMs.
Robbins conjectured the enumeration for odd order, but their number does not
appear to be round for even order.  No determinant or Pfaffian formula  is
known.

\item[$\bullet$] TSASMs - totally symmetric ASMs.  Their number does not appear
to be round and no determinant or Pfaffian formula is known.
\end{description}

Our results for these classes are as follows:

\begin{theorem}
The number of $n \times n$ ASMs is given by
$$A(n) = (-3)^{\binom{n}2}\prod_{i,j} \frac{3(j-i)+1}{j-i+n}.$$
The number of $2n+1 \times 2n+1$ vertically symmetric ASMs (VSASMs) is given by
$$A_\V(2n+1) =
(-3)^{n^2} \prod_{\substack{i,j \le 2n+1 \\ 2|j}}
\frac{3(j-i)+1}{j-i+2n+1}.$$
The number of $2n \times 2n$ half-turn symmetric ASMs (HTSASMs) is given by
$$\frac{A_\HT(2n)}{A(n)} =
(-3)^{\binom{n}2}\prod_{i,j} \frac{3(j-i)+2}{j-i+n}.$$
The number of $4n \times 4n$ quarter-turn symmetric ASMs (QTSASMs) is given by
$$A_\QT(4n) = A_\HT(2n)A(n)^2.$$
\eatline\label{th:main}
\end{theorem}

In the statement of \thm{th:main} and throughout this article, subscripts and
products range from 1 to $n$ unless otherwise specified.

\begin{theorem}  The 2- and 3-enumerations of ASMs and VSASMs
are given by
\begin{align*}
A(n;2) &= 2^{\binom{n}2} \\
A(n;3) &= \frac{3^{n^2-n}}{2^{n^2-n}}\prod_{\substack{i,j\\ 2\nmid j-i}}
    \frac{3(j-i)+1}{3(j-i)} \\
A_\V(2n+1;2) &= 2^{n^2-n} \\
A_\V(2n+1;3) &= \frac{3^{2n^2}}{2^{2n^2+n}}
    \prod_{\substack{i,j \le 2n+1 \\ 2\nmid i,2|j}} \frac{3(j-i)+1}{3(j-i)}.
\end{align*}
The 2-enumerations of even-sized HTSASMs and QTSASMs are given by 
\begin{align*}
A_\HT(2n;2,1) &= 2^{n^2} \prod_{\substack{i,j \\ 2\nmid j-i}}
    \frac{2(j-i)+1}{2(j-i)} \\
A_\QT(4n;2) &= (-1)^{\binom{n}2} 2^{2n^2-n} \prod_{i,j} \frac{4(j-i)+1}{j-i+n}.
\end{align*}
\eatline\label{th:3main}
\end{theorem}

\begin{fullfigure}{f:vhpasm}{The simplest VHPASM}
$$\begin{pmatrix} 0&0&0&+&0&0&0 \\ 0&+&0&-&0&+&0 \\
+&-&+&*&+&-&+ \\ 0&+&0&-&0&+&0 \\ 0&0&0&+&0&0&0
\end{pmatrix}$$ \eatline
\end{fullfigure}

\begin{fullfigure}{f:uasm}{A UASM}
$$\left(\begin{array}{ccc}
0&0&\rnode{a}{+} \\ 0&+&\rnode{b}{-} \\
+&-&\rnode{c}{0} \\ 0&0&\rnode{d}{+} \\
0&+&\rnode{e}{0} \\ 0&0&\rnode{f}{0}
\end{array}\right.\hspace{.3cm}.
\nccurve[angle=0,nodesepA=.2,nodesepB=.13,ncurv=1]{a}{b}
\nccurve[angle=0,nodesep=.2,ncurv=1]{c}{d}
\nccurve[angle=0,nodesep=.2,ncurv=1]{e}{f}$$
\eatline
\end{fullfigure}

We will also consider the following new types of ASMs:

\begin{description}
\item[$\bullet$] VHPASMs\footnote{Also known as $\beta$-ASMs, since their
boundary conditions are incompatible with VHS.} - vertically and horizontally
perverse ASMs.  A VHPASM has dimensions $4n+1 \times 4n+3$ for some integer
$n$.  It satisfies the alternating-sign condition and it has the same
symmetries as a VHSASM, except that the central entry ($*$) has the opposite
sign when read horizontally as when read vertically. The simplest VHPASM is
given in \fig{f:vhpasm}.

\item[$\bullet$] UASMs - $2n\times 2n$ ASMs with U-turn boundary on the right.
\fig{f:uasm} shows an example of a UASM.  As the example indicates, a UASM is
vertically just like an ASM. Horizontally the signs alternate if we read the
$2k-1$st row from left to right, and then continue to alternate if we read the
$2k$th row from right to left.  UASMs were first considered
by Tsuchiya \cite{Tsuchiya:determinant}.  They generalize VSASMs.

\item[$\bullet$] UUASMs\footnote{Also known as Unix-to-Unix ASMs.} - $2n\times
2n$ ASMs with U-turn boundary on the top and right.  UUASMs generalize both
VHSASMs and VHPASMs.

\item[$\bullet$] OSASMs - off-diagonally symmetric ASMs.  \Ie, DSASMs with
a null diagonal.

\item[$\bullet$] OOSASMs - off-diagonally, off-antidiagonally symmetric ASMs.
\Ie, DASASMs with null diagonals.

\item[$\bullet$] UOSASMs - off-diagonally symmetric UASMs.  They 
include TSASMs with null diagonals.
\end{description}

Our remaining results involve weighted enumerations that are more general than
just $x$-enumeration. We define the $y$-weight of a $2n \times 2n$ HTSASM to be
$y^k$ if the HTSASM has $k$ non-zero entries in the upper left quadrant.  This
yields the $(x,y)$-enumeration of HTSASMs, which is round when $y$ is $-1$ and
$x$ is 1 or 3.) We define the $x$-weight of a UASM or a UUASM be the number of
$-1$s, as before.  We define the $y$-weight of a UASM to be $y^k$ if $k$ of the
U-turns are oriented upward in the corresponding square ice state. We define
the $y$-weight of a UUASM the same way using the U-turns on the right, and
define the $z$-weight of a UUASM to be $z^k$ if $k$ of the U-turns on the top
are oriented to the right.  Thus we can consider the $(x,y)$-enumeration of
UASMs and the $(x,y,z)$-weight of UUASMs. As with UASMs, the $y$-weight of a
UOSASM is $y^k$ if  $k$ of the U-turns on the top are oriented to the right. 
By contrast, the $y$-weight of an OOSASM is $y^k$ if there are $2k$ more 1s
than $-1$s in the upper left quadrant.  Finally we index the generating
function of a given type of ASM by the length of one of its rows, counting the
length twice if the row takes a U-turn, and we include $x$-, $y$-, and
$z$-weight where applicable.  For example $A_\UO(8n;x,y)$ is the weighted
number of $4n \times 4n$ UOSASMs.

\begin{theorem}
There exist polynomials satisfying the equations
\begin{align*}
A(2n;x)                 &= 2A_\V(2n+1;x)\tA_\V(2n;x) \\
A(2n+1;x)               &= A_\V(2n+1;x)\tA_\V(2n+2;x) \\
A_\U(2n;x,y)            &= (y+1)^n A_\V(2n+1;x) \\
A_\HT(2n;x,\pm 1)       &= A(n;x) A^{(2)}_\HT(2n;x,\pm 1) \\
A_\UU(4n;x,y,z)         &= A_\V(2n+1;x) A^{(2)}_\UU(4n;x,y,z) \\
A_\OO(4n;x,y)           &= A_\O(2n;x) A^{(2)}_\OO(4n;x,y) \\
A_\QT(4n;x)             &= A^{(1)}_\QT(4n;x) A^{(2)}_\QT(4n;x) \\
A_\UO(8n;x,y)           &= A^{(1)}_\UO(8n;x) A^{(2)}_\UO(8n;x,y) \\
A^{(2)}_\HT(4n;x,1)     &= A^{(2)}_\UU(4n;x,1,1) \tA^{(2)}_\UU(4n;x) \\
A^{(2)}_\HT(4n+2;x,1)   &= 2A^{(2)}_\UU(4n;x,1,1) \tA^{(2)}_\UU(4n+4;x) \\
A^{(2)}_\HT(4n;x,-1)    &= (-x)^n A^{(1)}_\QT(4n;x)^2 \\
A^{(2)}_\OO(8n;x,-1)    &= (-x)^n A^{(1)}_\UO(8n;x,1)\tA^{(1)}_\UO(8n;x,1).
\end{align*}
\eatline\label{th:factor}
\end{theorem}

Many of the factorizations in \thm{th:factor} were conjectured experimentally
by David Robbins \cite{Robbins:symmetry}; the formula for $A_\U(2n;x,y)$ was
conjectured by Cohn and Propp \cite{CP:private}.

\begin{theorem}
The generating functions in \thm{th:factor} have the following special values.
\begin{align*}
A_\O(2n)              &= A_\V(2n+1) \\
A^{(2)}_\UU(4n;1,1,1) &= (-3)^{n^2}2^{2n}
    \prod_{\substack{i,j \le 2n+1 \\ 2|j}} \frac{3(j-i)+2}{j-i+2n+1}\\
A^{(2)}_\UU(4n;2,1,1) &= 2^{n(n+2)} \prod_{\substack{i,j \le 2n+1 \\ 2\nmid i,2|j}}
    \frac{2(j-i)+1}{2(j-i)} \\
A^{(2)}_\VHP(4n+2;1)  &= A_\V(2n+1) \\
A^{(1)}_\QT(4n)       &= A(n)^2 \\
A^{(1)}_\QT(4n;2)     &= (-1)^{\binom{n}2}2^{n(n-1)} \prod_{i,j} \frac{4(j-i)+1}{j-i+n} \\
A^{(1)}_\QT(4n;3)     &= 3^{\binom{n}2}A(n) \\
A^{(1)}_\UO(8n)       &= A_\V(2n+1)^2 \\
A^{(2)}_\UO(8n)       &= A_\UU(4n).
\end{align*}
\eatline\label{th:extra}
\end{theorem}

Other identities, for example that
$$A^{(2)}_\QT(4n)= A_\HT(2n),$$
are implied by combining Theorems~\ref{th:main}, \ref{th:3main},
\ref{th:factor}, and \ref{th:extra}, although such combinations do not always
reflect the logic of the proofs.

\section{Square ice}
\label{s:ice}

If $G$ is a tetravalent graph, an \emph{ice state} (also called a
\emph{six-vertex state}) of $G$ is an orientation of the edges such that two
edges enter and leave every tetravalent vertex.  In particular if $G$ is
locally a square grid, then the set of ice states is called \emph{square ice}
\cite{Lieb:exact}.  More generally $G$ may also have some univalent vertices,
which are called \emph{boundary}, and restrictions on the orientations of the
boundary edges are called \emph{boundary conditions}.

\begin{fullfigure}{f:grid}{Square ice with domain wall boundary}
\pspicture(-.5,-.5)(5,5)
\psline(0,1)(5,1) \psline(1,0)(1,5)
\psline(0,2)(5,2) \psline(2,0)(2,5)
\psline(0,3)(5,3) \psline(3,0)(3,5)
\psline(0,4)(5,4) \psline(4,0)(4,5)
\psline{->}(.61,1)(.62,1) \psline{->}(.61,2)(.62,2)
\psline{->}(.61,3)(.62,3) \psline{->}(.61,4)(.62,4)
\psline{->}(1,.39)(1,.38) \psline{->}(2,.39)(2,.38)
\psline{->}(3,.39)(3,.38) \psline{->}(4,.39)(4,.38)
\psline{->}(1,4.61)(1,4.62) \psline{->}(2,4.61)(2,4.62)
\psline{->}(3,4.61)(3,4.62) \psline{->}(4,4.61)(4,4.62)
\psline{->}(4.39,1)(4.38,1) \psline{->}(4.39,2)(4.38,2)
\psline{->}(4.39,3)(4.38,3) \psline{->}(4.39,4)(4.38,4)
\rput[r](-.2,1){$x_1$} \rput[t](1,-.2){$y_1$}
\rput[r](-.2,2){$x_2$} \rput[t](2,-.2){$y_2$}
\rput[r](-.2,3){$x_3$} \rput[t](3,-.2){$y_3$}
\rput[r](-.2,4){$x_4$} \rput[t](4,-.2){$y_4$}
\endpspicture
\end{fullfigure}

\begin{fullfigure}{f:replace}
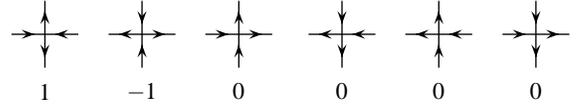
    {Replacing square ice with alternating-sign entries}
\verta{1}\hspace{.4cm}\vertb{-1}\hspace{.4cm}\vertc{0}\hspace{.4cm}
\vertd{0}\hspace{.4cm}\verte{0}\hspace{.4cm}\vertf{0}\eatline
\end{fullfigure}

For example, a finite square region of square ice can have \emph{domain wall
boundary}, defined as in at the sides and out at the top and bottom, as in
\fig{f:grid}.  These boundary conditions were first considered by Korepin
\cite{Korepin:norms,Izergin:finite,KBI:qism}.  A square ice state on this
region yields a matrix if we replace each vertex by a number according
\fig{f:replace}.  It is easy to check that this transformation is a bijection
between square ice with domain wall boundary and alternating-sign matrices
\cite{Kuperberg:asm,Kuperberg:eklp2}.

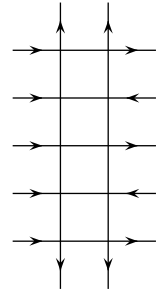
\begin{fullfigure}{f:vgrid}{Square ice with VS boundary}
\pspicture(0,0)(3.5,6)
\psline{->}(.61,1)(.62,1) \psline{->}(.61,2)(.62,2)
\psline{->}(.61,3)(.62,3) \psline{->}(.61,4)(.62,4)
\psline{->}(.61,5)(.62,5)
\psline{->}(1,.39)(1,.38) \psline{->}(2,.39)(2,.38)
\psline{->}(2.61,1)(2.62,1) \psline{->}(2.39,2)(2.38,2)
\psline{->}(2.61,3)(2.62,3) \psline{->}(2.39,4)(2.38,4)
\psline{->}(2.61,5)(2.62,5)
\psline{->}(1,5.61)(1,5.62) \psline{->}(2,5.61)(2,5.62)
\psline(0,1)(3,1) \psline(0,2)(3,2) \psline(0,3)(3,3)
\psline(0,4)(3,4) \psline(0,5)(3,5)
\psline(1,0)(1,6) \psline(2,0)(2,6)
\endpspicture
\end{fullfigure}

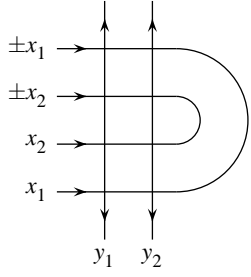
\begin{fullfigure}{f:htgrid}{Square ice with HTS boundary}
\pspicture(-.5,-.5)(4.5,5.5)
\psline(0,1)(2.5,1) \psline(0,2)(2.5,2)
\psline(0,3)(2.5,3) \psline(0,4)(2.5,4)
\psline(1,0)(1,5)   \psline(2,0)(2,5)
\psarc(2.5,2.5){0.5}{270}{90}
\psarc(2.5,2.5){1.5}{270}{90}
\psline{->}(.61,1)(.62,1) \psline{->}(1,.39)(1,.38)
\psline{->}(.61,2)(.62,2) \psline{->}(2,.39)(2,.38)
\psline{->}(.61,3)(.62,3) \psline{->}(1,4.61)(1,4.62)
\psline{->}(.61,4)(.62,4) \psline{->}(2,4.61)(2,4.62)
\rput[r](-.2,1){$x_1$} \rput[t](1,-.2){$y_1$}
\rput[r](-.2,2){$x_2$} \rput[t](2,-.2){$y_2$}
\rput[r](-.2,3){$\pm x_2$}
\rput[r](-.2,4){$\pm x_1$}
\endpspicture
\end{fullfigure}

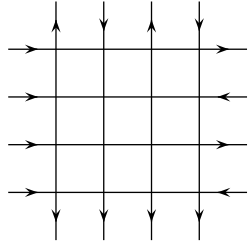
\begin{fullfigure}{f:vhsgrid}{Square ice with VHS boundary}
\pspicture(0,0)(5.5,5.5)
\psline(0,1)(5,1) \psline(1,0)(1,5)
\psline(0,2)(5,2) \psline(2,0)(2,5)
\psline(0,3)(5,3) \psline(3,0)(3,5)
\psline(0,4)(5,4) \psline(4,0)(4,5)
\psline{->}(.61,1)(.62,1) \psline{->}(.61,2)(.62,2)
\psline{->}(.61,3)(.62,3) \psline{->}(.61,4)(.62,4)
\psline{->}(1,.39)(1,.38) \psline{->}(2,.39)(2,.38)
\psline{->}(3,.39)(3,.38) \psline{->}(4,.39)(4,.38)
\psline{->}(4.39,1)(4.38,1) \psline{->}(1,4.61)(1,4.62)
\psline{->}(4.61,2)(4.62,2) \psline{->}(2,4.39)(2,4.38)
\psline{->}(4.39,3)(4.38,3) \psline{->}(3,4.61)(3,4.62)
\psline{->}(4.61,4)(4.62,4) \psline{->}(4,4.39)(4,4.38)
\endpspicture
\end{fullfigure}

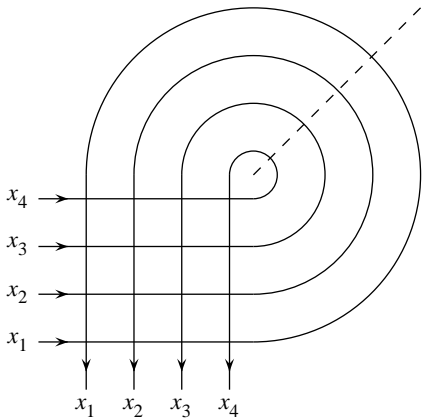
\begin{fullfigure}{f:qtgrid}{Square ice with QTS boundary}
\pspicture(-.5,-.5)(8,8)
\psline(0,1)(4.5,1) \psline(1,0)(1,4.5)
\psline(0,2)(4.5,2) \psline(2,0)(2,4.5)
\psline(0,3)(4.5,3) \psline(3,0)(3,4.5)
\psline(0,4)(4.5,4) \psline(4,0)(4,4.5)
\psarc(4.5,4.5){0.5}{270}{180}
\psarc(4.5,4.5){1.5}{270}{180}
\psarc(4.5,4.5){2.5}{270}{180}
\psarc(4.5,4.5){3.5}{270}{180}
\psline[linestyle=dashed](4.5,4.5)(8,8)
\psline{->}(.61,1)(.62,1) \psline{->}(.61,2)(.62,2)
\psline{->}(.61,3)(.62,3) \psline{->}(.61,4)(.62,4)
\psline{->}(1,.39)(1,.38) \psline{->}(2,.39)(2,.38)
\psline{->}(3,.39)(3,.38) \psline{->}(4,.39)(4,.38)
\rput[r](-.2,1){$x_1$} \rput[t](1,-.2){$x_1$}
\rput[r](-.2,2){$x_2$} \rput[t](2,-.2){$x_2$}
\rput[r](-.2,3){$x_3$} \rput[t](3,-.2){$x_3$}
\rput[r](-.2,4){$x_4$} \rput[t](4,-.2){$x_4$}
\endpspicture
\end{fullfigure}

There are also easy bijections from ice states of the graphs in
Figures~\ref{f:vgrid}, \ref{f:htgrid}, \ref{f:vhsgrid}, and \ref{f:qtgrid} to
the sets of VSASMs, VHSASMs, even HTSASMs, and even QTSASMs.  (The labels in
these figures will be used later.) The dashed line in the QTSASM graph means
that the orientation of an edge reverses as it crosses the line. The HTSASM and
QTSASM graphs are obtained by quotienting the unrestricted ASM graph by the
symmetry.  The median of a $2n+1 \times 2n+1$ VSASM is always the same, so we
can delete it and consider the alternating-sign patterns on the left half. The
deleted median then produces the alternating boundary in \fig{f:vgrid}. 
Likewise we can quarter a VHSASM by deleting both medians, which produces two
alternating sides.

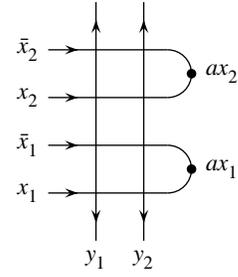
\begin{fullfigure}{f:ugrid}{Square ice with U boundary}
\pspicture(-.5,-.5)(3.6,5)
\psline{->}(.61,1)(.62,1) \psline{->}(.61,2)(.62,2)
\psline{->}(.61,3)(.62,3) \psline{->}(.61,4)(.62,4)
\psline{->}(1,.39)(1,.38) \psline{->}(2,.39)(2,.38)
\psline{->}(1,4.61)(1,4.62) \psline{->}(2,4.61)(2,4.62)
\psline(0,1)(2.5,1) \psline(0,2)(2.5,2)
\psline(0,3)(2.5,3) \psline(0,4)(2.5,4)
\psline(1,0)(1,5) \psline(2,0)(2,5)
\psarc(2.5,1.5){.5}{270}{90} \psarc(2.5,3.5){.5}{270}{90}
\rput[r](-.2,1){$x_1$} \rput[t](1,-.2){$y_1$}
\rput[r](-.2,2){$\bx_1$} \rput[t](2,-.2){$y_2$}
\rput[r](-.2,3){$x_2$}
\rput[r](-.2,4){$\bx_2$}
\qdisk(3,1.5){.1} \rput[l](3.3,1.5){$ax_1$}
\qdisk(3,3.5){.1} \rput[l](3.3,3.5){$ax_2$}
\endpspicture
\end{fullfigure}

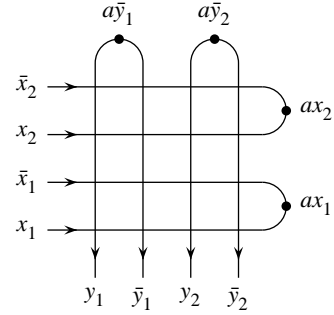
\begin{fullfigure}{f:uugrid}{Square ice with UU boundary}
\pspicture(-.5,-.5)(5.6,5.6)
\psline(0,1)(4.5,1) \psline(1,0)(1,4.5)
\psline(0,2)(4.5,2) \psline(2,0)(2,4.5)
\psline(0,3)(4.5,3) \psline(3,0)(3,4.5)
\psline(0,4)(4.5,4) \psline(4,0)(4,4.5)
\psarc(1.5,4.5){.5}{0}{180}  \psarc(3.5,4.5){.5}{0}{180}
\psarc(4.5,1.5){.5}{270}{90} \psarc(4.5,3.5){.5}{270}{90}
\psline{->}(.61,1)(.62,1) \psline{->}(.61,2)(.62,2)
\psline{->}(.61,3)(.62,3) \psline{->}(.61,4)(.62,4)
\psline{->}(1,.39)(1,.38) \psline{->}(2,.39)(2,.38)
\psline{->}(3,.39)(3,.38) \psline{->}(4,.39)(4,.38)
\rput[r](-.2,1){$x_1$}   \rput[t](1,-.2){$y_1$}
\rput[r](-.2,2){$\bx_1$} \rput[t](2,-.2){$\by_1$}
\rput[r](-.2,3){$x_2$}   \rput[t](3,-.2){$y_2$}
\rput[r](-.2,4){$\bx_2$} \rput[t](4,-.2){$\by_2$}
\qdisk(1.5,5){.1} \rput[b](1.5,5.3){$a\by_1$}
\qdisk(3.5,5){.1} \rput[b](3.5,5.3){$a\by_2$}
\qdisk(5,1.5){.1} \rput[l](5.3,1.5){$ax_1$}
\qdisk(5,3.5){.1} \rput[l](5.3,3.5){$ax_2$}
\endpspicture
\end{fullfigure}

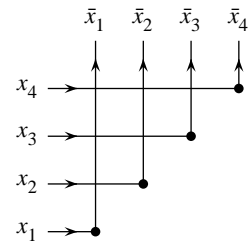
\begin{fullfigure}{f:ogrid}{Square ice with OS boundary}
\pspicture(-.5,.5)(5,6)
\psline(0,1)(1,1)(1,5) \psline(0,2)(2,2)(2,5)
\psline(0,3)(3,3)(3,5) \psline(0,4)(4,4)(4,5)
\psline{->}(.61,1)(.62,1) \psline{->}(.61,2)(.62,2)
\psline{->}(.61,3)(.62,3) \psline{->}(.61,4)(.62,4)
\psline{->}(1,4.61)(1,4.62) \psline{->}(2,4.61)(2,4.62)
\psline{->}(3,4.61)(3,4.62) \psline{->}(4,4.61)(4,4.62)
\rput[r](-.2,1){$x_1$} \rput[b](1,5.2){$\bx_1$}
\rput[r](-.2,2){$x_2$} \rput[b](2,5.2){$\bx_2$}
\rput[r](-.2,3){$x_3$} \rput[b](3,5.2){$\bx_3$}
\rput[r](-.2,4){$x_4$} \rput[b](4,5.2){$\bx_4$}
\qdisk(1,1){.1} \qdisk(2,2){.1} \qdisk(3,3){.1} \qdisk(4,4){.1}
\endpspicture
\end{fullfigure}

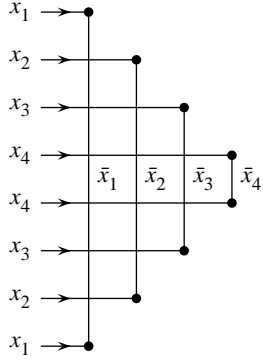
\begin{fullfigure}{f:oogrid}{Square ice with OOS boundary}
\pspicture(-.5,-.5)(4.5,8)
\psline(0,1)(1,1)(1,8)(0,8) \psline(0,2)(2,2)(2,7)(0,7)
\psline(0,3)(3,3)(3,6)(0,6) \psline(0,4)(4,4)(4,5)(0,5)
\psline{->}(.61,1)(.62,1) \psline{->}(.61,2)(.62,2)
\psline{->}(.61,3)(.62,3) \psline{->}(.61,4)(.62,4)
\psline{->}(.61,5)(.62,5) \psline{->}(.61,6)(.62,6)
\psline{->}(.61,7)(.62,7) \psline{->}(.61,8)(.62,8)
\rput[r](-.2,1){$x_1$} \rput[r](-.2,2){$x_2$}
\rput[r](-.2,3){$x_3$} \rput[r](-.2,4){$x_4$}
\rput[r](-.2,5){$x_4$} \rput[r](-.2,6){$x_3$}
\rput[r](-.2,7){$x_2$} \rput[r](-.2,8){$x_1$}
\rput[l](1.2,4.5){$\bx_1$} \rput[l](2.2,4.5){$\bx_2$}
\rput[l](3.2,4.5){$\bx_3$} \rput[l](4.2,4.5){$\bx_4$}
\qdisk(1,1){.1} \qdisk(2,2){.1} \qdisk(3,3){.1} \qdisk(4,4){.1}
\qdisk(4,5){.1} \qdisk(3,6){.1} \qdisk(2,7){.1} \qdisk(1,8){.1}
\endpspicture
\end{fullfigure}

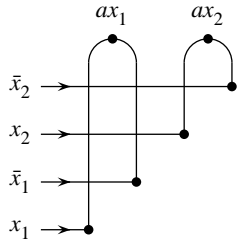
\begin{fullfigure}{f:uogrid}{Square ice with UOS boundary}
\pspicture(-.5,-.5)(4.5,6)
\psline(0,1)(1,1) \psline(0,2)(2,2) \psline(0,3)(3,3) \psline(0,4)(4,4)
\psline(1,1)(1,4.5) \psline(2,2)(2,4.5)
\psline(3,3)(3,4.5) \psline(4,4)(4,4.5)
\psarc(1.5,4.5){.5}{0}{180} \psarc(3.5,4.5){.5}{0}{180}
\psline{->}(.61,1)(.62,1) \psline{->}(.61,2)(.62,2)
\psline{->}(.61,3)(.62,3) \psline{->}(.61,4)(.62,4)
\rput[r](-.2,1){$x_1$} \rput[r](-.2,2){$\bx_1$}
\rput[r](-.2,3){$x_2$} \rput[r](-.2,4){$\bx_2$}
\qdisk(1,1){.1} \qdisk(2,2){.1} \qdisk(3,3){.1} \qdisk(4,4){.1}
\qdisk(1.5,5){.1} \rput[b](1.5,5.3){$ax_1$}
\qdisk(3.5,5){.1} \rput[b](3.5,5.3){$ax_2$}
\endpspicture
\end{fullfigure}

Finally the square ice grids corresponding to UASMs, UUASMs, OSASMs, OOSASMs,
and UOSASMs are shown in Figures~\ref{f:ugrid}, \ref{f:uugrid}, \ref{f:ogrid},
\ref{f:oogrid}, and \ref{f:uogrid}. The last three grids have right-angled
divalent vertices; we require the orientations of a square ice state to either
be both in or both out at these vertices.  In contrast at the U-turn
vertices one edge must point in and one must point out.

\section{Local concerns}
\label{s:local}

Throughout the article we assume the following abbreviations:
\begin{align*}
\bx &= x^{-1} \\
\sigma(x) &= x - \bx \\
\alpha(x) &= \sigma(ax)\sigma(a\bx).
\end{align*}
(As we discuss below, $a$ is a global parameter that need not
appear as an explicit argument of $\alpha$.)

We will consider a class of multiplicative weights for symmetric ASMs. By a
multiplicative weight we mean that the weight of some object is the product of
the weights of its parts.  In statistical mechanics, multiplicative weights are
called \emph{Boltzmann weights}, and the total weight of all objects is called
a \emph{partition function}. \fig{f:weights} shows the weights that we will use
for the six possible states of a vertex.  The figure also shows the weights for
U-turns and corners that are labelled with a dot; bare edges and curves have
the trivial weight 1. The vertex weights are called an \emph{$R$-matrix} and
the U-turn and corner weights are called \emph{$K$-matrices}.  Vertex and
U-turn weights depend on a parameter $x$ (the \emph{spectral parameter}) which
may be different for different vertices or U-turns, so we will label sites by
the value of $x$.  The weights also depend on three parameters $a$, $b$, and
$c$ which will be the same for all elements of any single square ice grid, so
these parameters do not appear as labels.

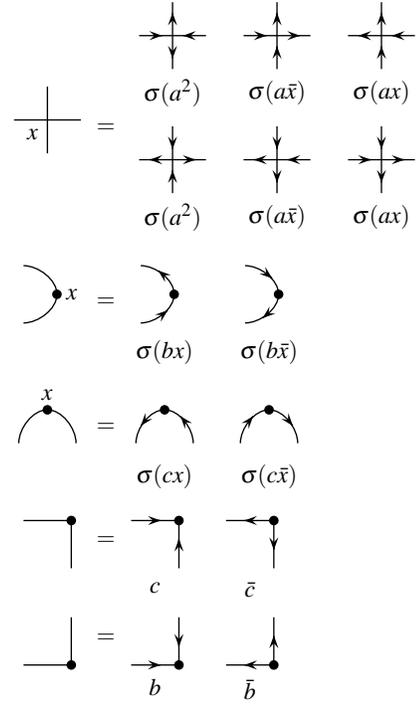
\begin{fullfigure}{f:weights}{Weights for vertices and U-turns}
\begin{align*}
\pspicture[.6](-.7,-1.7)(.7,.7)
\psline(0,0)(0,.7)    \psline(0,0)(.7,0)
\psline(-0,0)(-.7,0)  \psline(0,-0)(0,-.7)
\rput[tr](-.2,-.2){$x$}\endpspicture\;\;&=\;\;
\begin{array}{c@{\hspace{.5cm}}c@{\hspace{.5cm}}c}
\verta{\sigma(a^2)}&\vertc{\sigma(a\bx)}&\verte{\sigma(ax)} \\
\vertb{\sigma(a^2)}&\vertd{\sigma(a\bx)}&\vertf{\sigma(ax)}
\end{array} \\
\pspicture[.6](-.7,-1.7)(.7,.7)
\pccurve[angleA=90,angleB=0](.2,0)(-.5,.6)
\pccurve[angleA=270,angleB=0](.2,0)(-.5,-.6)\qdisk(.2,0){.1}
\rput[l](.4,0){$x$}\endpspicture\;\;&=\;\;
\pspicture[.6](-.7,-1.7)(.7,.7)
\pccurve[angleA=90,angleB=0](.2,0)(-.5,.6)\mto
\pccurve[angleA=270,angleB=0](.2,0)(-.5,-.6)\mfro\qdisk(.2,0){.1}
\rput(0,-1.2){$\sigma(bx)$}\endpspicture\hspace{.5cm}
\pspicture[.6](-.7,-1.7)(.7,.7)
\pccurve[angleA=90,angleB=0](.2,0)(-.5,.6)\mfro
\pccurve[angleA=270,angleB=0](.2,0)(-.5,-.6)\mto\qdisk(.2,0){.1}
\rput(0,-1.2){$\sigma(b\bx)$}\endpspicture \\
\pspicture[.6](-.7,-1.7)(.7,.7)
\pccurve[angleA=0,angleB=90](0,.2)(.6,-.5)
\pccurve[angleA=180,angleB=90](0,.2)(-.6,-.5)\qdisk(0,.2){.1}
\rput[b](0,.4){$x$}\endpspicture\;\;&=\;\;
\pspicture[.6](-.7,-1.7)(.7,.7)
\pccurve[angleA=0,angleB=90](0,.2)(.6,-.5)\mfro
\pccurve[angleA=180,angleB=90](0,.2)(-.6,-.5)\mto\qdisk(0,.2){.1}
\rput(0,-1.2){$\sigma(cx)$}\endpspicture\hspace{.5cm}
\pspicture[.6](-.7,-1.7)(.7,.7)
\pccurve[angleA=0,angleB=90](0,.2)(.6,-.5)\mto
\pccurve[angleA=180,angleB=90](0,.2)(-.6,-.5)\mfro\qdisk(0,.2){.1}
\rput(0,-1.2){$\sigma(c\bx)$}\endpspicture \\
\pspicture[.6](0,-.6)(1.2,1.2)
\psline(1,0)(1,1)(0,1)\qdisk(1,1){.1}\endpspicture\;\;&=\;\;
\pspicture[.6](0,-.6)(1.2,1.2)
\psline(1,0)(1,1)(0,1)\qdisk(1,1){.1}\rput[t](.5,-.3){$c$}
\psline{->}(1,.61)(1,.62)
\psline{->}(.61,1)(.62,1)
\endpspicture
\hspace{.5cm}
\pspicture[.6](0,-.6)(1.2,1.2)
\psline(1,0)(1,1)(0,1)\qdisk(1,1){.1}\rput[t](.5,-.3){$\bc$}
\psline{->}(1,.39)(1,.38)
\psline{->}(.39,1)(.38,1)
\endpspicture \\
\pspicture[.6](0,-.6)(1.2,1.2)
\psline(0,0)(1,0)(1,1)\qdisk(1,0){.1}\endpspicture\;\;&=\;\;
\pspicture[.6](0,-.6)(1.2,1.2)
\psline(0,0)(1,0)(1,1)\qdisk(1,0){.1}\rput[t](.5,-.3){$b$}
\psline{->}(1,.39)(1,.38)
\psline{->}(.61,0)(.62,0)
\endpspicture
\hspace{.5cm}
\pspicture[.6](0,-.6)(1.2,1.2)
\psline(0,0)(1,0)(1,1)\qdisk(1,0){.1}\rput[t](.5,-.3){$\bb$}
\psline{->}(1,.61)(1,.62)
\psline{->}(.39,0)(.38,0)
\endpspicture
\end{align*}\eatline
\end{fullfigure}

We will use a graph with labelled vertices as a notation for its corresponding
partition function.  If the graph has unoriented boundary edges, then the
partition function is also interpreted as a function of the orientations of the
edges.  On the other hand, our definitions imply that we sum over the
orientations of internal edges.  For example, the graph
$$
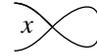
\pspicture(-1,-.5)(1,.5)
\pccurve[angleA=0,angleB=225](-.8,-.5)(0,0)
\pccurve[angleA=0,angleB=135](-.8,.5)(0,0)
\psbezier(-.8,.5)(-.5,.5)(-.3,.3)(0,0)
\psbezier(0,0)(.4,.4)(1,.7)(1,0)
\psbezier(0,0)(.4,-.4)(1,-.7)(1,0)
\rput[r](-.4,0){$x$}
\endpspicture
$$
denotes the following function on the set of four orientations of the
boundary:
$$
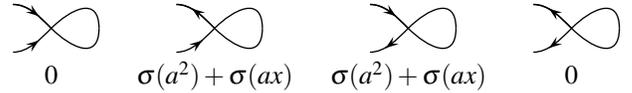
\begin{array}{c@{\hspace{.5cm}}c@{\hspace{.5cm}}c@{\hspace{.5cm}}c}
\pspicture(-1,-.5)(1,.5)
\pccurve[angleA=0,angleB=225](-.8,-.5)(0,0)\mto
\pccurve[angleA=0,angleB=135](-.8,.5)(0,0)\mto
\psbezier(-.8,.5)(-.5,.5)(-.3,.3)(0,0)
\psbezier(0,0)(.4,.4)(1,.7)(1,0)
\psbezier(0,0)(.4,-.4)(1,-.7)(1,0)
\endpspicture &
\pspicture(-1,-.5)(1,.5)
\pccurve[angleA=0,angleB=225](-.8,-.5)(0,0)\mto
\pccurve[angleA=0,angleB=135](-.8,.5)(0,0)\mfro
\psbezier(-.8,.5)(-.5,.5)(-.3,.3)(0,0)
\psbezier(0,0)(.4,.4)(1,.7)(1,0)
\psbezier(0,0)(.4,-.4)(1,-.7)(1,0)
\endpspicture &
\pspicture(-1,-.5)(1,.5)
\pccurve[angleA=0,angleB=225](-.8,-.5)(0,0)\mfro
\pccurve[angleA=0,angleB=135](-.8,.5)(0,0)\mto
\psbezier(-.8,.5)(-.5,.5)(-.3,.3)(0,0)
\psbezier(0,0)(.4,.4)(1,.7)(1,0)
\psbezier(0,0)(.4,-.4)(1,-.7)(1,0)
\endpspicture &
\pspicture(-1,-.5)(1,.5)
\pccurve[angleA=0,angleB=225](-.8,-.5)(0,0)\mfro
\pccurve[angleA=0,angleB=135](-.8,.5)(0,0)\mfro
\psbezier(-.8,.5)(-.5,.5)(-.3,.3)(0,0)
\psbezier(0,0)(.4,.4)(1,.7)(1,0)
\psbezier(0,0)(.4,-.4)(1,-.7)(1,0)
\endpspicture \\
0 & \sigma(a^2)+\sigma(ax) & \sigma(a^2)+\sigma(ax) & 0
\end{array}
$$

In this notation a vertex is not quite invariant under rotation by
90 degrees, so the meaning of a label depends on the quadrant in which
it appears.  The following relation holds:
$$\pspicture[.45](-.7,-.7)(.7,.7)
\psline(0,0)(0,.7)\psline(0,0)(.7,0)\psline(-0,0)(-.7,0)\psline(0,-0)(0,-.7)
\rput[tr](-.2,-.2){$x$}\endpspicture\;\;=\;\;
\pspicture[.45](-.7,-.7)(.7,.7)
\psline(0,0)(0,.7)\psline(0,0)(.7,0)\psline(-0,0)(-.7,0)\psline(0,-0)(0,-.7)
\rput[tl](.2,-.2){$\bx$}\endpspicture$$

As a further abbreviation, if we label two lines of a graph that cross at an
unlabelled vertex, the spectral parameter is set to their ratio:
$$\pspicture[.6](-.7,-1.5)(.7,.7)
\psline(0,0)(0,.7)\psline(0,0)(.7,0)\psline(-0,0)(-.7,0)\psline(0,-0)(0,-.7)
\rput[r](-1,0){$x$}\rput[t](0,-1){$y$}
\endpspicture\;\;=\;\;
\pspicture[.6](-.7,-1.5)(.7,.7)
\psline(0,0)(0,.7)\psline(0,0)(.7,0)\psline(-0,0)(-.7,0)\psline(0,-0)(0,-.7)
\rput[tr](-.2,-.2){$x\by$}\endpspicture$$

The labelled graphs in Figures~\ref{f:grid}, \ref{f:vgrid}, \ref{f:htgrid},
\ref{f:qtgrid}, \ref{f:ugrid}, \ref{f:uugrid}, \ref{f:ogrid}, \ref{f:oogrid},
and \ref{f:uogrid} then represent the
partition functions
$$\begin{array}{cccc}
Z(n;\vx,\vy) & Z^\pm_\HT(n;\vx,\vy) & Z_\U(n;\vx,\vy) & Z_\UU(n;\vx,\vy) \\
Z_\QT(n;\vx) & Z_\O(n;\vx) & Z_\OO(n;\vx) & Z_\UO(n;\vx)
\end{array}$$
Here the vectors $\vx$ and $\vy$ have length $n$ when both are present,
and otherwise $\vx$ has length $2n$.  In the HT and OO cases
there is a single extra parameter taken from the set $\{+,-\}$; if it is
$-$ then the spectral parameters in the upper half of the grid are
negated.  (Note that the index $n$ is not defined in the same way
as for the enumerators such as $A_\HT(2n)$.)

The key property of the $R$-matrix is that it satisfies the
Yang-Baxter equation:

\begin{lemma}[Yang-Baxter equation]  If $xyz = \ba$, then
$$\pspicture[.42](-1.8,-1)(1.1,2)
\pnode(0,0){b1}\pnode(0,1){b2}\pnode(-.866,.5){b3}
\pnode([angle=255,nodesep=.8]b1){a1}\pnode([angle=345,nodesep=.8]b1){a2}
\pnode([angle= 15,nodesep=.8]b2){a3}\pnode([angle=105,nodesep=.8]b2){a4}
\pnode([angle=135,nodesep=.8]b3){a5}\pnode([angle=225,nodesep=.8]b3){a6}
\pscurve(a1)(b1)(b2)(a4)\pscurve(a3)(b2)(b3)(a6)\pscurve(a5)(b3)(b1)(a2)
\rput[tl](.1,-.3){$x$}\rput[bl](.1,1.3){$y$}\rput[r](-1.2,.5){$z$}
\endpspicture = \pspicture[.42](-1.1,-1)(1.8,2)
\pnode(0,0){b1}\pnode(0,1){b2}\pnode(.866,.5){b3}
\pnode([angle=285,nodesep=.8]b1){a1}\pnode([angle=195,nodesep=.8]b1){a2}
\pnode([angle=165,nodesep=.8]b2){a3}\pnode([angle= 75,nodesep=.8]b2){a4}
\pnode([angle= 45,nodesep=.8]b3){a5}\pnode([angle=315,nodesep=.8]b3){a6}
\pscurve(a1)(b1)(b2)(a4)\pscurve(a3)(b2)(b3)(a6)\pscurve(a5)(b3)(b1)(a2)
\rput[tr](-.1,-.3){$y$}\rput[br](-.1,1.3){$x$}\rput[l](1.2,.5){$z$}
\endpspicture\hspace{.1cm}.$$
\label{l:ybe}\eatline
\end{lemma}

As usual the Yang-Baxter equation appears to be a massive coincidence. In our
previous review of the Yang-Baxter equation \cite{Kuperberg:asm}, the
$R$-matrix was normalized to have a particular symmetry:  It was the matrix of
an invariant tensor over the 2-dimensional representation of the quantum group
$U_q(\sl(2))$, with $q$ related to our present parameter $a$.  This
symmetry reduced the coincidence in the equation to a single numerical
equality.  The spectral parameters were chosen to satisfy the equality.  Here
we normalize the $R$-matrix to reveal combinatorial symmetry rather than
symmetry from quantum algebra.

\begin{proof} Taken literally, the equation consists of 64 numerical
equalities, because there are 64 ways to orient the six boundary edges on each
side.   However, both sides are zero unless three edges point in and three
point out.  This leaves 20 non-zero equations.  The equation also has three
kinds of symmetry:  The right side is the left side rotated by 180 degrees, all
arrows may be reversed, and both sides may be rotated by 120 degrees if the
variables $x$, $y$, and $z$ are cyclically permuted. By the three symmetries, 8
of the non-zero equations are tautological, and the other 12 are all
equivalent.  One of the 12 non-trivial equations is
$$\pspicture[.42](-1.8,-1)(1.1,2)
\pnode(0,0){b1}\pnode(0,1){b2}\pnode(-.866,.5){b3}
\pnode([angle=255,nodesep=.8]b1){a1}\pnode([angle=345,nodesep=.8]b1){a2}
\pnode([angle= 15,nodesep=.8]b2){a3}\pnode([angle=105,nodesep=.8]b2){a4}
\pnode([angle=135,nodesep=.8]b3){a5}\pnode([angle=225,nodesep=.8]b3){a6}
\ncline[linestyle=none]{a1}{b1}\mto\ncline[linestyle=none]{a2}{b1}\mto
\ncline[linestyle=none]{a3}{b2}\mfro\ncline[linestyle=none]{a4}{b2}\mfro
\ncline[linestyle=none]{a5}{b3}\mto\ncline[linestyle=none]{a6}{b3}\mfro
\pscurve(a1)(b1)(b2)(a4)\pscurve(a3)(b2)(b3)(a6)\pscurve(a5)(b3)(b1)(a2)
\rput[tl](.1,-.3){$x$}\rput[bl](.1,1.3){$y$}\rput[r](-1.2,.5){$z$}
\endpspicture = \pspicture[.42](-1.1,-1)(1.8,2)
\pnode(0,0){b1}\pnode(0,1){b2}\pnode(.866,.5){b3}
\pnode([angle=285,nodesep=.8]b1){a1}\pnode([angle=195,nodesep=.8]b1){a2}
\pnode([angle=165,nodesep=.8]b2){a3}\pnode([angle= 75,nodesep=.8]b2){a4}
\pnode([angle= 45,nodesep=.8]b3){a5}\pnode([angle=315,nodesep=.8]b3){a6}
\ncline[linestyle=none]{a1}{b1}\mto\ncline[linestyle=none]{a2}{b1}\mfro
\ncline[linestyle=none]{a3}{b2}\mto\ncline[linestyle=none]{a4}{b2}\mfro
\ncline[linestyle=none]{a5}{b3}\mfro\ncline[linestyle=none]{a6}{b3}\mto
\pscurve(a1)(b1)(b2)(a4)\pscurve(a3)(b2)(b3)(a6)\pscurve(a5)(b3)(b1)(a2)
\rput[tr](-.1,-.3){$y$}\rput[br](-.1,1.3){$x$}\rput[l](1.2,.5){$z$}
\endpspicture\hspace{.1cm}.$$
In algebraic form, the equation is
$$\sigma(a\by)\sigma(a^2)\sigma(a\bx) =
\sigma(az)\sigma(a^2)^2 + \sigma(ax)\sigma(ay)\sigma(a^2).$$
Cancelling a factor of $\sigma(a^2)$, expanding, and cancelling terms
yields
$$a^2\bx\by + \ba^2xy = a^3z -z\ba - \bz a + \ba^3\bz + a^2xy  +
\ba^2\bx\by,$$
which is implied by the condition $xyz=\ba$.
\end{proof}

We will need the reflection equation
\cite{Cherednik:factorizing,Sklyanin:boundary,VG:boundary}, an analogue of the
Yang-Baxter equation that relates a $K$-matrix to the $R$-matrix.

\begin{lemma}[Reflection equation] If $st=ay$ and $s\bt=ax$, then
$$
\pspicture[.45](-.4,-1.6)(.8,1.6)
\pccurve[angleA=270,angleB=135](-.4,1.2)(0,.6)
\pccurve[angleA=270,angleB=45](.4,1.2)(0,.6)
\pccurve[angleA=315,angleB=90](0,.6)(.4,0)
\pccurve[angleA=225,angleB=90](0,.6)(-.4,0)
\pccurve[angleA=270,angleB=135](-.4,0)(0,-.6)
\pccurve[angleA=270,angleB=45](.4,0)(0,-.6)
\pccurve[angleA=315,angleB=90](0,-.6)(.4,-1.2)
\pccurve[angleA=225,angleB=90](0,-.6)(-.4,-1.2)
\psline(.4, 1.2)(.4, 1.6)\psline(-.4, 1.2)(-.4, 1.6)
\psline(.4,-1.2)(.4,-1.6)\psline(-.4,-1.2)(-.4,-1.6)
\qdisk(.4,0){.1}\rput[l](.7,0){$s$}
\qdisk(.4,-1.2){.1}\rput[l](.7,-1.2){$t$}
\rput[b](0,1){$x$}\rput[t](0,-1){$y$}
\endpspicture\;\;=\;\;
\pspicture[.45](-.4,-1.6)(.8,1.6)
\pccurve[angleA=270,angleB=135](-.4,1.2)(0,.6)
\pccurve[angleA=270,angleB=45](.4,1.2)(0,.6)
\pccurve[angleA=315,angleB=90](0,.6)(.4,0)
\pccurve[angleA=225,angleB=90](0,.6)(-.4,0)
\pccurve[angleA=270,angleB=135](-.4,0)(0,-.6)
\pccurve[angleA=270,angleB=45](.4,0)(0,-.6)
\pccurve[angleA=315,angleB=90](0,-.6)(.4,-1.2)
\pccurve[angleA=225,angleB=90](0,-.6)(-.4,-1.2)
\psline(.4, 1.2)(.4, 1.6)\psline(.4,-1.2)(.4,-1.6)
\psline(-.4, 1.2)(-.4, 1.6)\psline(-.4,-1.2)(-.4,-1.6)
\qdisk(.4,0){.1}\rput[l](.7,0){$s$}
\qdisk(.4,1.2){.1}\rput[l](.7,1.2){$t$}
\rput[b](0,1){$y$}\rput[t](0,-1){$x$}
\endpspicture
$$
\label{l:reu}\eatline
\end{lemma}
\begin{proof}  The argument is similar to that for Lemma~\ref{l:ybe}.
Both sides are zero unless two boundary edges point in and two
point out.   There is a symmetry exchanging the two sides given
by reflecting through a horizontal line and simultaneously reversing
all arrows.  (Note that the weights of a U-turn are not invariant
under reflection alone.)  Under this symmetry 4 of the 6 non-zero
equations are tautological, and the other 2 are equivalent.
One of these is:
$$
\pspicture[.45](-.4,-1.8)(.8,1.8)
\pccurve[angleA=270,angleB=135](-.4,1.2)(0,.6)
\pccurve[angleA=270,angleB=45](.4,1.2)(0,.6)
\pccurve[angleA=315,angleB=90](0,.6)(.4,0)
\pccurve[angleA=225,angleB=90](0,.6)(-.4,0)
\pccurve[angleA=270,angleB=135](-.4,0)(0,-.6)
\pccurve[angleA=270,angleB=45](.4,0)(0,-.6)
\pccurve[angleA=315,angleB=90](0,-.6)(.4,-1.2)
\pccurve[angleA=225,angleB=90](0,-.6)(-.4,-1.2)
\pcline(-.4, 1.2)(-.4, 1.8)\mto\pcline(.4, 1.2)(.4, 1.8)\mfro
\pcline(-.4,-1.2)(-.4,-1.8)\mto\pcline(.4,-1.2)(.4,-1.8)\mfro
\qdisk(.4,0){.1}\rput[l](.7,0){$s$}
\qdisk(.4,-1.2){.1}\rput[l](.7,-1.2){$t$}
\rput[b](0,1){$x$}\rput[t](0,-1){$y$}
\endpspicture\;\;=\;\;
\pspicture[.45](-.4,-1.6)(.8,1.6)
\pccurve[angleA=270,angleB=135](-.4,1.2)(0,.6)
\pccurve[angleA=270,angleB=45](.4,1.2)(0,.6)
\pccurve[angleA=315,angleB=90](0,.6)(.4,0)
\pccurve[angleA=225,angleB=90](0,.6)(-.4,0)
\pccurve[angleA=270,angleB=135](-.4,0)(0,-.6)
\pccurve[angleA=270,angleB=45](.4,0)(0,-.6)
\pccurve[angleA=315,angleB=90](0,-.6)(.4,-1.2)
\pccurve[angleA=225,angleB=90](0,-.6)(-.4,-1.2)
\pcline(-.4, 1.2)(-.4, 1.8)\mto\pcline(.4, 1.2)(.4, 1.8)\mfro
\pcline(-.4,-1.2)(-.4,-1.8)\mto\pcline(.4,-1.2)(.4,-1.8)\mfro
\qdisk(.4,0){.1}\rput[l](.7,0){$s$}
\qdisk(.4,1.2){.1}\rput[l](.7,1.2){$t$}
\rput[b](0,1){$y$}\rput[t](0,-1){$x$}
\endpspicture
$$
Algebraically, the equation reads:
\begin{align*}
\sigma(bt)\sigma(a^2)(\sigma(ay)&\sigma(b\bs) + \sigma(ax)\sigma(bs)) = \\
&\sigma(b\bt)\sigma(a^2)(\sigma(ax)\sigma(b\bs) + \sigma(ay)\sigma(bs)).
\end{align*}
All terms of the equation match or cancel when $st=ay$ and $s\bt=ax$.
\end{proof}

The corner $K$-matrix also satisfies the reflection equation
\cite{VG:boundary}.

\begin{lemma} For any $x$ and $y$,
$$
\pspicture[.45](-.5,0)(2,2)
\qdisk(2,1){.1} \qdisk(1,2){.1}
\psline(1,1)(2,1)(2,0) \psline(1,2)(1,0)
\psbezier(-.5,2)(0,2)(.5,1)(1,1)
\psbezier(-.5,1)(0,1)(.5,2)(1,2)
\rput[r](-.2,1.5){$x$}\rput[rt](.8,.8){$y$}
\endpspicture\;\;=\;\;
\pspicture[.45](0,-.5)(2,2)
\qdisk(1,2){.1} \qdisk(2,1){.1}
\psline(1,1)(1,2)(0,2) \psline(2,1)(0,1)
\psbezier(2,-.5)(2,0)(1,.5)(1,1)
\psbezier(1,-.5)(1,0)(2,.5)(2,1)
\rput[t](1.5,-.2){$x$}\rput[rt](.8,.8){$y$}
\endpspicture\;\;.
$$
\label{l:rec}
\end{lemma}
\begin{proof} Diagonal reflection exchanges the two sides.  Both
sides are zero if an odd number of boundary edges point inward.  If two
boundary edges point in and the other two point out, then arrow reversal is
also a symmetry, because one corner must have inward arrows and the other
outward arrows.  These facts together imply that all cases of the  equation are
null or tautological.
\end{proof}

Finally we will need an equation that, loosely speaking, inverts a U-turn:

\begin{lemma}[Fish equation] For any $a$ and $x$,
$$
\pspicture[.4](-1,-.5)(1.7,.5)
\pccurve[angleA=0,angleB=225](-.8,-.5)(0,0)
\pccurve[angleA=0,angleB=135](-.8,.5)(0,0)
\psbezier(-.8,.5)(-.5,.5)(-.3,.3)(0,0)
\psbezier(0,0)(.4,.4)(1,.7)(1,0)\psbezier(0,0)(.4,-.4)(1,-.7)(1,0)
\rput[r](-.5,0){$\ba x^2$}\qdisk(1,0){.1}\rput[l](1.3,0){$ax$}
\endpspicture\;\;=\;\;\sigma(a^2x^2)
\pspicture[.4](-1,-.5)(1.1,.5)
\psline(-.5,.5)(0,.5)\psline(-.5,-.5)(0,-.5)\psarc(0,0){.5}{270}{90}
\qdisk(.5,0){.1}\rput[l](.8,0){$a\bx$}
\endpspicture
$$
\label{l:fish}
\eatline
\end{lemma}

The proof is elementary.

\section{Determinants}
\label{s:det}

In this section we will establish determinant
and Pfaffian formulas for the partition functions defined in
Sections~\ref{s:ice} and \ref{s:local}.  Recall that 
the Pfaffian of an antisymmetric $2n \times 2n$
matrix $A$ is defined as
$$\Pf A \stackrel{\mathrm{def}}{=}
\sum_{\pi \in X} (-1)^\pi \prod_i
A_{\pi(2i-1),\pi(2i)},$$
where $X \subset S_{2n}$ has one representative in each coset of the wreath
product $S_2 \wr S_n$. (Thus $X$ admits a bijection with the set of perfect
matchings of $\{1,\ldots,2n\}$.)  Recall also that 
$$\det A = (\Pf A)^2.$$

\begin{theorem} Let 
\begin{align*}
M(n;\vx,\vy)_{i,j} =&\  \frac1{\alpha(x_i\by_j)} \\
M^\pm_\HT(n;\vx,\vy)_{i,j} =&\  \frac1{\sigma(a\bx_iy_j)} \pm 
\frac1{\sigma(ax_i\by_j)} \\
M_\U(n;\vx,\vy)_{i,j} =&\  \frac1{\alpha(x_i\by_j)}-\frac1{\alpha(x_iy_j)} \\
M_\UU(n;\vx,\vy)_{i,j} =&\ 
\frac{\sigma(b\by_j)\sigma(cx_i)}{\sigma(ax_i\by_j)} -
\frac{\sigma(b\by_j)\sigma(c\bx_i)}{\sigma(a\bx_i\by_j)} \\
&\ - \frac{\sigma(by_j)\sigma(cx_i)}{\sigma(ax_iy_j)} +
\frac{\sigma(by_j)\sigma(c\bx_i)}{\sigma(a\bx_iy_j)} \\
M^{(k)}_\QT(n;\vx)_{i,j} =&\  \frac{\sigma(\bx_i^kx_j^k)}{\alpha(\bx_ix_j)} \\
M_\O(n;\vx)_{i,j} =& \frac{\sigma(\bx_ix_j)}{\alpha(x_ix_j)} \\
M_\OO(n;\vx)_{i,j} =& \sigma(\bx_i x_j)\biggl(\frac{c^2}{\sigma(ax_ix_j)}
    + \frac{b^2}{\sigma(a\bx_i\bx_j)}\biggr) \\
M^{(1)}_\UO(n;\vx)_{i,j} =& \sigma(\bx_ix_j)\sigma(x_ix_j)
    \biggl(\frac1{\alpha(x_ix_j)}-\frac1{\alpha(\bx_ix_j)}\biggr) \\
M^{(2)}_\UO(n;\vx)_{i,j} =& \sigma(\bx_ix_j)\sigma(x_ix_j)
    \biggl(\frac{\sigma(cx_i)\sigma(cx_j)}{\sigma(ax_ix_j)}
    - \frac{\sigma(cx_i)\sigma(c\bx_j)}{\sigma(ax_i\bx_j)} \\
    &\ - \frac{\sigma(c\bx_i)\sigma(cx_j)}{\sigma(a\bx_ix_j)}
    + \frac{\sigma(c\bx_i)\sigma(c\bx_j)}{\sigma(a\bx_i\bx_j)}\biggr).
\end{align*}
Then
\begin{align*}
Z(n;\vx,\vy) =&\  \frac{\sigma(a^2)^n\prod_{i,j} \alpha(x_i\by_j)}
    {\prod_{i<j}\sigma(\bx_ix_j)\sigma(y_i\by_j)}(\det M) \\
Z^\pm_\HT(n;\vx,\vy) =&\  \frac{\sigma(a^2)^n \prod_{i,j}\alpha(x_i\by_j)^2}
    {\prod_{i<j}\sigma(\bx_ix_j)^2\sigma(y_i\by_j)^2}
    (\det M)(\det M^\pm_\HT) \\
Z_\U(n;\vx,\vy) =&\  \frac{\sigma(a^2)^n\prod_i \sigma(b\by_i)\sigma(a^2x_i^2)
    \prod_{i,j}\alpha(x_i\by_j)\alpha(x_iy_j)}
    {\prod_{i<j}\sigma(\bx_ix_j)\sigma(y_i\by_j)
    \prod_{i\le j}\sigma(\bx_i\bx_j)\sigma(y_iy_j)} \\
    &\ \cdot (\det M_\U)
\end{align*}
\begin{align*}
Z_\UU(n;\vx,\vy) =&\ \frac{\sigma(a^2)^n\prod_i \sigma(a^2x_i^2)
    \sigma(a^2\by_i^2)\prod_{i,j}\alpha(x_i\by_j)^2\alpha(x_iy_j)^2}
    {\prod_{i<j}\sigma(\bx_ix_j)^2\sigma(y_i\by_j)^2
    \prod_{i\le j}\sigma(\bx_i\bx_j)^2\sigma(y_iy_j)^2} \\
    &\ \cdot (\det M_\U)(\det M_\UU) \\
Z_\QT(n;\vx) =&\  \frac{\sigma(a^2)^n\sigma(a)^{3n}\prod_{i<j \le 2n}
    \alpha(\bx_ix_j)^2}{\prod_{i<j \le 2n} \sigma(\bx_ix_j)^2} \\
    &\ \cdot (\Pf M^{(1)}_\QT)(\Pf M^{(2)}_\QT) \\
Z_\O(n;\vx) =&\ \frac{\sigma(a^2)^n\prod_{i<j \le 2n}\alpha(x_ix_j)}
    {\prod_{i<j \le 2n} \sigma(\bx_ix_j)} (\Pf M_\O)\\
Z_\OO(n;\vx) =&\ \frac{\sigma(a^2)^n\prod_{i<j \le 2n}\alpha(x_ix_j)^2}
    {\prod_{i<j \le 2n} \sigma(\bx_ix_j)^2} (\Pf M_\O)(\Pf M_\OO) \\
Z_\UO(n;\vx) =&\ \frac{b^{2n}\sigma(a^2)^n\sigma(a)^{2n}\prod_{i \le 2n}
    \sigma(a^2x_i^2)}{\prod_{i < j \le 2n}\sigma(\bx_i x_j)^2
    \prod_{i\le j\le 2n}\sigma(x_i x_j)^2} \\
    &\ \cdot \!\!\!\prod_{i<j\le 2n}\!\!\!\!\alpha(\bx_i x_j)^2
    \alpha(x_i x_j)^2 (\Pf M^{(1)}_\UO)(\Pf M^{(2)}_\UO)
\end{align*}\eatline
\label{th:z}
\end{theorem}

We call the first four partition functions the \emph{determinant partition
functions} and the other four the \emph{Pfaffian partition functions}.

\begin{remark} The partition function $Z_\U(n;\vx,\vy)$, the Tsuchiya
determinant, is nearly invariant if $\vx$ is exchanged with $\vy$.  Similarly
$Z_\O(n;\vx)$ is nearly invariant if each $x_i$ is replaced with $\bx_i$. We
have no direct explanation for these symmetries.  Note that the first symmetry
is less apparent in Tsuchiya's matrix $M$ \cite[eq.
(42)]{Tsuchiya:determinant}, which has an asymmetric factor
$$F_{ij} = \frac{\sinh(\zeta_- + \lambda_j)}{\sinh(\lambda_j + \omega_i)}
+ \frac{\sinh(\zeta_- - \lambda_j)}{\sinh(\lambda_j - \omega_i)}.$$
In this expression $\omega_i$, $\lambda_j$, and $\zeta_-$ are
obtained from $x_i$, $y_j$, and $b$ by reparameterization.
If we factor this expression,
$$F_{ij} = \frac{\sinh(2\lambda_j)\sinh(\zeta_- - \omega_i)}
    {\sinh(\lambda_j + \omega_i)\sinh(\lambda_j - \omega_i)}.$$
we can then pull the asymmetric factors out of the determinant since
they each depend on only one of the two indices $i$ and $j$.
This also explains why the $K$-matrix parameter $\zeta_-$ or $b$
need not appear in the matrix $M_U$.
\end{remark}

The proof of \thm{th:z} uses recurrence relations that determine both sides.
The relations are expressed in Lemmas~\ref{l:sym}, \ref{l:xysym},
\ref{l:spec}, and \ref{l:poly}. Indeed, the first three of these lemmas are
obvious for the right-hand sides of \thm{th:z}; only Lemma~\ref{l:poly} needs
to be argued for both sides.

\begin{lemma}[Baxter, Sklyanin] Each of the partition functions in
\thm{th:z} is symmetric in the coordinates of $\vx$.  Each determinant
partition function is symmetric in the coordinates of $\vy$.  The partition
functions $Z_\U(n;\vx,\vy)$ and $Z_\UU(n;\vx,\vy)$ gain a factor of
$\sigma(a^2\bx_i^2)/\sigma(a^2x_i^2)$ if $x_i$ is replaced by $\bx_i$
for a single $i$.  Similarly $Z_\UU(n;\vx,\vy)$ gains
$\sigma(a^2y_i^2)/\sigma(a^2\by_i^2)$ under $y_i \mapsto \by_i$
and $Z_\UO(n;\vx)$ gains $\sigma(a^2\bx_i^2)/\sigma(a^2x_i^2)$
under $x_i \mapsto \bx_i$.
\label{l:sym}
\end{lemma}
\begin{proof} Invariance of $Z(n;\vx,\vy)$ is an illustrative case. We exchange
$x_i$ with $x_{i+1}$ for any $i \le n-1$ by crossing the corresponding lines at
the left side.  If the spectral parameter of the crossing is $z = \ba
x_i\bx_{i+1}$, we can move it to the right side using the Yang-Baxter equation
(Lemma~\ref{l:ybe}) and then remove it:
\begin{align*}
& \sigma(az)
\pspicture[.45](-1.75,.25)(6,2.75)
\rput[r](-.7,1){$x_i$}
\rput[r](-.7,2){$x_{i+1}$}
\psline{->}(.11,1)(.12,1) \psline{->}(.11,2)(.12,2)
\psline(-.5,1)(2,1) \psline(-.5,2)(2,2)
\psline(3,1)(4.5,1) \psline(3,2)(4.5,2)
\rput(2.5,1.5){\large $\ldots$}
\psline(.5,.5)(.5,2.5) \psline(1.5,.5)(1.5,2.5) \psline(3.5,.5)(3.5,2.5)
\psline{->}(3.89,1)(3.88,1) \psline{->}(3.89,2)(3.88,2)
\endpspicture \\
&\qquad\qquad= \pspicture[.45](-2.25,.25)(5.5,2.75)
\rput[r](-1.2,2){$x_i$}
\rput[r](-1.2,1){$x_{i+1}$}
\psbezier(-.5,2)(0,2)(.5,1)(1,1)
\psbezier(-.5,1)(0,1)(.5,2)(1,2)
\rput[r](-.1,1.5){$z$}
\psline(1,1)(2.5,1) \psline(1,2)(2.5,2)
\rput(3,1.5){\large $\ldots$}
\psline(3.5,1)(5,1) \psline(3.5,2)(5,2)
\psline(1,.5)(1,2.5) \psline(2,.5)(2,2.5) \psline(4,.5)(4,2.5)
\psline{->}(-1,1)(-.5,1) \psline{->}(-1,2)(-.5,2)
\psline{->}(4.39,1)(4.38,1) \psline{->}(4.39,2)(4.38,2)
\endpspicture \\
&\qquad\qquad= \pspicture[.45](-1.75,.25)(5.75,2.75)
\rput[r](-.7,2){$x_i$}
\rput[r](-.7,1){$x_{i+1}$}
\psbezier(.5,2)(1,2)(1.5,1)(2,1)
\psbezier(.5,1)(1,1)(1.5,2)(2,2)
\rput[r](.9,1.5){$z$}
\psline(-.5,1)(.5,1) \psline(-.5,2)(.5,2)
\psline(2,1)(2.5,1) \psline(2,2)(2.5,2)
\rput(3,1.5){\large $\ldots$}
\psline(3.5,1)(5,1) \psline(3.5,2)(5,2)
\psline(.5,.5)(.5,2.5) \psline(2,.5)(2,2.5) \psline(4,.5)(4,2.5)
\psline{->}(.11,1)(.12,1) \psline{->}(.11,2)(.12,2)
\psline{->}(4.39,1)(4.38,1) \psline{->}(4.39,2)(4.38,2)
\endpspicture \\
&\qquad\qquad= \pspicture[.45](-1.75,.25)(6,2.75)
\rput[r](-.7,2){$x_i$}
\rput[r](-.7,1){$x_{i+1}$}
\psline(-.5,1)(2,1) \psline(-.5,2)(2,2)
\psline(3,1)(3.5,1) \psline(3,2)(3.5,2)
\psbezier(3.5,2)(4,2)(4.5,1)(5,1)
\psbezier(3.5,1)(4,1)(4.5,2)(5,2)
\rput[r](3.9,1.5){$z$}
\rput(2.5,1.5){\large $\ldots$}
\psline{<-}(5,1)(5.5,1) \psline{<-}(5,2)(5.5,2)
\psline(.5,.5)(.5,2.5) \psline(1.5,.5)(1.5,2.5) \psline(3.5,.5)(3.5,2.5)
\psline{->}(.11,1)(.12,1) \psline{->}(.11,2)(.12,2)
\endpspicture \\
&\qquad\qquad= \sigma(az)\pspicture[.45](-1.75,.25)(5,2.75)
\rput[r](-.7,2){$x_i$}
\rput[r](-.7,1){$x_{i+1}$}
\psline{->}(.11,1)(.12,1) \psline{->}(.11,2)(.12,2)
\psline(-.5,1)(2,1) \psline(-.5,2)(2,2)
\psline(3,1)(4.5,1) \psline(3,2)(4.5,2)
\rput(2.5,1.5){\large $\ldots$}
\psline(.5,.5)(.5,2.5) \psline(1.5,.5)(1.5,2.5) \psline(3.5,.5)(3.5,2.5)
\psline{->}(3.89,1)(3.88,1) \psline{->}(3.89,2)(3.88,2)
\endpspicture\;\;.
\end{align*}
The argument for symmetry in $\vx$ is exactly the same for all of the square
ice grids without U-turns.  If the grid has diagonal boundary with corner
vertices, we can bounce the crossing off of it using Lemma~\ref{l:rec}. 

If the grid has U-turn boundary on the right, we exchange $x_i$ with $x_{i+1}$
by  crossing the $\bx_i$ line over the two lines above it.  We let the spectral
parameters of these two crossings be $z = \ba \bx_i x_{i+1}$ and $w = \ba \bx_i
\bx_{i+1}$.  We move both crossings to the right using the Yang-Baxter
equation, then we bounce them off of the U-turns using the reflection
equation (Lemma~\ref{l:reu}):
\begin{align*}
&\pspicture[.48](-2.5,0)(7,4.75)
\psline{->}(-.39,1)(-.38,1) \psline{->}(-.39,2)(-.38,2)
\psline{->}(-.39,3)(-.38,3) \psline{->}(-.39,4)(-.38,4)
\psline(-1,1)(.5,1) \psline(1.5,1)(5.5,1)
\psline(-1,2)(.5,2) \psline(1.5,2)(4,2)
\psline(-1,3)(.5,3) \psline(1.5,3)(2.5,3)
\psline(-1,4)(.5,4) \psline(1.5,4)(2.5,4)
\rput(1,1.5){\large $\ldots$}
\rput(1,3.5){\large $\ldots$}
\psbezier(4,2)(4.5,2)(5,3)(5.5,3) \psbezier(4,3)(4.5,3)(5,2)(5.5,2)
\psbezier(2.5,3)(3,3)(3.5,4)(4,4) \psbezier(2.5,4)(3,4)(3.5,3)(4,3)
\rput[r](4.4,2.5){$w$} \rput[r](2.9,3.5){$z$}
\psline(4,4)(5.5,4)
\psline(0,.5)(0,4.5) \psline(2,.5)(2,4.5)
\psarc(5.5,1.5){.5}{270}{90} \psarc(5.5,3.5){.5}{270}{90}
\rput[r](-1.2,1){$x_i$}
\rput[r](-1.2,2){$x_{i+1}$}
\rput[r](-1.2,3){$\bx_{i+1}$}
\rput[r](-1.2,4){$\bx_i$}
\qdisk(6,1.5){.1} \rput[l](6.3,1.5){$ax_i$}
\qdisk(6,3.5){.1} \rput[l](6.3,3.5){$ax_{i+1}$}
\endpspicture \\
&\qquad\quad= \pspicture[.48](-2.5,0)(7,4.75)
\psline{->}(-.39,1)(-.38,1) \psline{->}(-.39,2)(-.38,2)
\psline{->}(-.39,3)(-.38,3) \psline{->}(-.39,4)(-.38,4)
\psline(-1,1)(.5,1) \psline(1.5,1)(2.5,1)
\psline(-1,2)(.5,2) \psline(1.5,2)(2.5,2)
\psline(-1,3)(.5,3) \psline(1.5,3)(4,3)
\psline(-1,4)(.5,4) \psline(1.5,4)(5.5,4)
\rput(1,1.5){\large $\ldots$}
\rput(1,3.5){\large $\ldots$}
\psbezier(4,2)(4.5,2)(5,3)(5.5,3) \psbezier(4,3)(4.5,3)(5,2)(5.5,2)
\psbezier(2.5,2)(3,2)(3.5,1)(4,1) \psbezier(2.5,1)(3,1)(3.5,2)(4,2)
\rput[r](4.4,2.5){$w$} \rput[r](2.9,1.5){$z$}
\psline(4,1)(5.5,1)
\psline(0,.5)(0,4.5) \psline(2,.5)(2,4.5)
\psarc(5.5,1.5){.5}{270}{90} \psarc(5.5,3.5){.5}{270}{90}
\rput[r](-1.2,1){$x_i$}
\rput[r](-1.2,2){$x_{i+1}$}
\rput[r](-1.2,3){$\bx_{i+1}$}
\rput[r](-1.2,4){$\bx_i$}
\qdisk(6,1.5){.1} \rput[l](6.3,1.5){$ax_{i+1}$}
\qdisk(6,3.5){.1} \rput[l](6.3,3.5){$ax_i$}
\endpspicture\;\;.
\end{align*}
Also if the grid has a U-turn on the right, we establish covariance under $x_i
\mapsto \bx_i$ by switching the lines with these two labels and eating the
crossing using the fish equation (Lemma~\ref{l:fish}).

The same arguments establish symmetry in the coordinates of $\vy$.
All of the arguments used in combination establish the claimed properties
of $Z_\UO(n;\vx)$.
\end{proof}

\begin{lemma} The partition function $Z_\HT(n;\vx,\vy)$ gains
a factor of $(\pm 1)^n$ if $x_i$ and $y_i$ are replaced by $\bx_i$ and
$\by_i$ for all $i$ simultaneously.  Similarly $Z_{\OO,\pm}(n;\vx)$
is invariant under $x_i \mapsto \bx_i$ and
$b \mapsto c \mapsto b$.
\label{l:xysym}
\end{lemma}
\begin{proof} In both cases, the symmetry is effected by reflecting the square
ice grid or the alternating-sign matrices through a horizontal line.
\end{proof}

For a vector $\vx = (x_1,\ldots,x_n)$, let $\vx' = (x_2,\ldots,x_n)$.

\begin{lemma}  If $x_1 = ay_1$, then
\begin{align*}
\frac{Z(n;\vx,\vy)}{Z(n-1;\vx',\vy')} =&\ \sigma(a^2)\prod_{2 \le i}
    \sigma(a\bx_iy_1)\sigma(a\bx_1y_i)  \\
\frac{Z^\pm_\HT(n;\vx,\vy)}{Z^\pm_\HT(n-1;\vx',\vy')} =&\
    \pm \sigma(a^2)^2\prod_{2\le i}\sigma(a\bx_iy_1)^2\sigma(a\bx_1y_i)^2 \\
\frac{Z_\U(n;\vx,\vy)}{Z_\U(n-1;\vx',\vy')} =&\
    \sigma(a^2) \sigma(a^2x_1^2)\sigma(b\by_1) \\[-2ex]
    &\ \!\!\cdot \prod_{2 \le i} \sigma(a\bx_iy_1)\sigma(a\bx_1y_i)
    \sigma(ax_iy_1)\sigma(a\bx_1\by_i) \\
\frac{Z_\UU(n;\vx,\vy)}{Z_\UU(n-1;\vx',\vy')} 
    =&\ \sigma(a^2)^2 \sigma(a^2x_1^2)\sigma(a^2\by_1^2)
    \sigma(b\by_1) \sigma(cx_1) \\[-2ex]
    &\ \cdot \prod_{2 \le i} \sigma(a\bx_iy_1)^2\sigma(a\bx_1y_i)^2 \\
    &\ \cdot \prod_{2 \le i} \sigma(ax_iy_1)^2\sigma(a\bx_1\by_i)^2.
\end{align*}
If $x_2 = ax_1$, then
$$
\frac{Z_\QT(n;\vx)}{Z_\QT(n-1;\vx'')} = \sigma(a)^2\sigma(a^2)^2
    \!\!\prod_{3 \le i \le 2n}\!\!\!\! \sigma(ax_i\bx_1)^2\sigma(a\bx_ix_2)^2.
$$
If $x_2 = \ba\bx_1$, then
\begin{align*}
\frac{Z_\O(n;\vx)}{Z_\O(n-1;\vx'')} =&\ \sigma(a^2)
    \!\!\prod_{3 \le i \le 2n}\!\!\!\!\sigma(a\bx_1\bx_i)\sigma(a\bx_2\bx_i) \\
\frac{Z_\OO(n;\vx)}{Z_\OO(n-1;\vx'')}
    =&\ c^2\sigma(a^2)^2\!\!\prod_{3 \le i \le 2n}\!\!\!\!
    \sigma(a\bx_1\bx_i)^2\sigma(a\bx_2\bx_i)^2 \\
\frac{Z_\UO(n;\vx)}{Z_\UO(n-1;\vx'')} =&\ b^2\sigma(a)^2\sigma(a^2)^2
    \sigma(a^2x_1^2)\sigma(a^2x_2^2) \\[-2ex]
    & \cdot \sigma(cx_1)\sigma(cx_2)\!\!\prod_{3 \le i \le 2n}\!\!\!\!
    \sigma(a\bx_1\bx_i)^2 \sigma(a\bx_2\bx_i)^2 \\
    & \cdot \!\!\prod_{3 \le i \le 2n}\!\!\!\!
    \sigma(a\bx_1x_i)^2\sigma(a\bx_2x_i)^2.
\end{align*}\eatline
\label{l:spec}
\end{lemma}
\begin{proof} This lemma is clearer in the alternating-sign matrix model than
it is in the square ice model.  The partition function $Z(n;\vx,\vy)$ is a sum
over $n \times n$ alternating-sign matrices in which each entry of the matrix
has a multiplicative weight. When $y_1 = ax_1$, the weight of a 0 in the
southwest corner is 0. Consequently this corner is forced to be 1 and the left
column and bottom row are forced to be 0, as in \fig{f:forced}.  The sum
reduces to  one over $(n-1) \times (n-1)$ ASMs.  The only discrepancy between
$Z(n;\vx,\vy)|_{y_1 = ax_1}$ and $Z(n-1;\vx',\vy')$  is the weights of
the forced entries, which the lemma lists as factors.

\begin{fullfigure}{f:forced}{ASM entries forced by $y_1 = ax_1$}
$$\left(\begin{array}{c|cccc}
0&0&0&+&0 \\ 0&+&0&-&+ \\ 0&0&0&+&0 \\ 0&0&+&0&0 \\ \hline +&0&0&0&0
\end{array}\right)$$
\eatline
\end{fullfigure}

The argument in the other determinant cases is identical. The argument in the
Pfaffian cases is only slightly different:  All QTSASMs have zeroes in the
corners, and the specialization $x_2 = ax_1$ instead forces a 1 next to each
corner and zeroes the first two rows and columns from each edge.  Likewise the
specialization $x_2 = \ba\bx_1$ forces a 1 next to each corner of an OSASM or
an OOSASM and a 1 in the third row entry bottom of a UOSASM, and several rows
and columns of zeroes in each of these cases.
\end{proof}

Define the \emph{width} of a Laurent polynomial to be the
difference in degree between the leading and trailing terms. (For example,
$q^3-q^{-2}$ has width 5.)

\begin{lemma} Both sides of each equation of \thm{th:z} are Laurent polynomials
in each coordinate of $\vx$ (and $\vy$ in the determinant cases) and their
widths in $x_1$ ($y_1$ in the determinant cases) are as given in
\tab{t:widths}.
\label{l:poly}
\end{lemma}

To conclude the proof of \thm{th:z}, we claim that Lemmas~\ref{l:sym},
\ref{l:xysym}, and \ref{l:spec} inductively determine both sides by Lagrange
interpolation.  (To begin the induction each partition function is set to 1
when $n=0$.)   If a Laurent polynomial of width $w$ has prespecified leading
and trailing exponents, it is determined by $w+1$ of its values. Each of our
partition functions is a centered Laurent polynomial in $x_1$ (in the Pfaffian
cases) or $y_1$ (in the determinant cases).  Moreover each is either an even
function or an odd function.  Thus we only need $w+1$
specializations, where $w$ is the width in $x_1^2$ (or $y_1^2$).  

These widths are summarized in \tab{t:widths}.  To compute them, observe that
each 0 entry in the bottom row of an ASM contributes 1 to the width.  In the
UASM and UUASM cases, it is the bottom two rows, and the U-turn itself
contributes 1 to the width as well.  In the QTSASM case, the corner entries
always have weight $\sigma(a)$ and do not contribute to the width.
Lemmas~\ref{l:sym} and \ref{l:spec} together provide many specializations which
are listed in \tab{t:widths}.  Note that Lemma~\ref{l:sym} implies that
$\sigma(a^2x_1^2)$ divides $Z_\UU(n;\vx,\vy)$ $Z_\UO(n;\vx,vy)$, which provides
an extra specialization in these two cases.  In conclusion, it is easy to check
that there are enough specializations to match the widths. 

\begin{fulltable}{t:widths}{Widths and specializations of partition functions.}
\begin{tabular}{c@{\hspace{.5cm}}c@{\hspace{.5cm}}l}
Function & Width & Specializations\\ \hline
$Z(n;\vx,\vy)$         & $n-1$  & $y_1 = ax_i$ \\
$Z^\pm_\HT(n;\vx,\vy)$ & $2n-1$ & $y_1 = a^{\pm 1} x_i$ \\
$Z_\U(n;\vx,\vy)$      & $2n-1$ & $y_1 = ax_i^{\pm 1}$ \\
$Z_\UU(n;\vx,\vy)$     & $4n$   & $y_1 = a^{\pm 1}x_i^{\pm 1}$, $a$ \\
$Z_\QT(n;\vx)$         & $4n-3$ & $x_1 = a^{\pm 1} x_i$ \\
$Z_\O(n;\vx)$          & $2n-2$ & $x_1 = a\bx_i$ \\
$Z_\OO(n;\vx)$         & $4n-3$ & $x_1 = a^{\pm 1}\bx_i$ \\
$Z_\UO(n;\vx)$         & $8n-4$ & $x_1 = a^{\pm 1}x_i^{\pm 1}$, $\ba$
\end{tabular}\eatline
\end{fulltable}

\begin{remark} The formulas in \thm{th:z} are even more special than
Lemmas~\ref{l:sym} through \ref{l:poly} suggest.  Among the evidence for this,
the recurrence relations still hold with only slight modifications if all
spectral parameters in the QT, UU, and UO grids are multiplied by an extra
parameter $z$.  Similarly the spectral parameters in the top halves of the HT
and OO grids may be multiplied by an arbitrary $z$ instead of by $\pm 1$. 
However, we were not able to generalize \thm{th:z} to include this parameter.

Lemma~\ref{l:spec} reveals another subtlety, namely that
$$\frac{Z^+_\HT(n;\vx,\vy)}{Z^+_\HT(n-1;\vx',\vy')}
    = \left(\frac{Z(n;\vx,\vy)}{Z(n-1;\vx',\vy')}\right)^2$$
at every specialization $y_i = a^{\pm 1} x_j$.  Since this coincidence
holds for enough specializations to determine $Z^+_\HT(n;\vx,\vy)$
entirely, one might suppose that
$$Z^+_\HT(n;\vx,\vy) = Z(n;\vx,\vy)^2.$$
But then $Z^+_\HT(n;\vx,\vy)$ would be an even function of $y_1$, while in
reality it is an odd function.  The other symmetry classes involving half-turn
rotation have similar behavior.
\end{remark}

\section{Factor exhaustion}
\label{s:factor}

In this section we derive several round determinants and Pfaffians depending on
two and three parameters.  We will later identify special cases of the
determinants and Pfaffians with those appearing in \thm{th:z}, and they will
specialize further to establish the enumerations in Theorems~\ref{th:main},
\ref{th:3main}, and \ref{th:extra}.

Since the formulas in this section may seem complicated, we recommend verifying
that they are round without worrying about their exact form in the first
reading.  For this purpose we give a more precise definition of roundness that
also applies to polynomials.  A term $R_n$ in a sequence of rational
polynomials depending on one or more variables is \emph{round} if it is a ratio
of products of constants, monomials, and differences of two monic terms.  All
exponents and constant factors should grow polynomially in $n$ or be
independent of $n$. For example, $n!3^n(q+p^n)$ is round and Gaussian binomial
coefficients are round. Note that a round polynomial in a single variable must
be a product of cyclotomic polynomials, which is part of the motivation for the
term ``round''.  Roundness is preserved when a variable is set to 1 or to a
product of other variables.  In a later reading one can verify the explicit
formulas. This is a tedious but elementary computation, because  all round
expressions involved have an explicit and regular form.  As a warmup the reader
can verify that the expressions for $A(n)$ in Theorems~\ref{th:asm} and
\ref{th:main} coincide.

We begin with the classic Cauchy double alternant and a Pfaffian
generalization found independently by Stembridge
and by Laksov, Lascoux, and Thorup
\cite{Knuth:pfaffians,Stembridge:paths,LLT:giambelli,Sundquist:pfaffian}.

\begin{theorem}[Cauchy, S. L. L. T.] Let
\begin{align*}
C_1(\vx,\vy)_{i,j} = \frac1{x_i+y_j} \\
C_2(\vx,\vy)_{i,j} = \frac1{x_i+y_j} - \frac1{1+x_iy_j}
\end{align*}
For $1 \le i,j \le 2n$, let 
\begin{align*}
C_3(\vx)_{i,j} &= \frac{x_j-x_i}{x_i + x_j} \\
C_4(\vx)_{i,j} &= \frac{x_j-x_i}{1-x_ix_j}.
\end{align*}
Then
\begin{align*}
\det C_1 =&\ \frac{\prod_{i<j} (x_j-x_i)(y_j-y_i)}{\prod_{i,j} (x_i+y_j)} \\
\det C_2 =&\ \frac{\prod_{i<j} (1-x_ix_j)(1-y_iy_j)(x_j-x_i)(y_j-y_i)}{\prod_{i,j} (x_i+y_j)(1+x_iy_j)}\\
          &\ \cdot \prod_i (1+x_i)(1+y_i) \\
\Pf C_3  =&\ \prod_{i<j\le 2n} \frac{x_i-x_j}{x_i + x_j} \\
\Pf C_4  =&\ \prod_{i<j\le 2n} \frac{x_i-x_j}{1-x_ix_j}
\end{align*}\eatline\label{th:cauchy}
\end{theorem}
\begin{proof} Our proof is by the factor exhaustion method
\cite{Krattenthaler:advanced}. The determinant $\det C_1$ is divisible by $x_j-x_i$
because when $x_i = x_j$, two rows of $C_1$ are proportional.  Likewise it is
also divisible by $y_j-y_i$. At the same time, the polynomial
$$\prod_{i,j} (x_i+y_j)(\det C_1)$$
has degree $n^2-n$, so it has no room for other non-constant factors. This
determines $\det C_1$ up to a constant, which can be found inductively
by setting $x_1 = -y_1$.

The determinant $\det C_2$ is argued the same way. The Pfaffians $\Pf C_3$ and
$\Pf C_4$ are also argued the same way; here the constant factor can be found
by setting $x_1 = \bx_2$.
\end{proof}

Next we evaluate four determinants in the variables $p$ and $q$.  We use two
more functions similar to $\sigma$ and $\alpha$ from Section~\ref{s:det}:
$$\gamma(q) = q^{1/2} - q^{-1/2} \qquad \tau(q) = q^{1/2} + q^{-1/2}.$$

\begin{theorem} Let
\begin{align*}
T_1(p,q)_{i,j} &= \frac{\gamma(q^{n+j-i})}{\gamma(p^{n+j-i})} \\
T_2(p,q)_{i,j} &= \frac{\tau(q^{j-i})}{\tau(p^{j-i})} \\
T_3(p,q)_{i,j} &= \frac{\gamma(q^{n+j+i})}{\gamma(p^{n+j+i})} -
\frac{\gamma(q^{n+j-i})}{\gamma(p^{n+j-i})} \\
T_4(p,q)_{i,j} &= \frac{\tau(q^{j+i})}{\tau(p^{j+i})} -
\frac{\tau(q^{j-i})}{\tau(p^{j-i})}
\end{align*}
Then
\begin{align*}
\det T_1 &= \frac{\prod_{i \ne j} \gamma(p^{j-i})\prod_{i,j} \gamma(qp^{j-i})}
    {\prod_{i,j} \gamma(p^{n+j-i})}\\
\det T_2 &= (-1)^{\binom{n}2}
    \frac{2^n\prod_{\substack{i \ne j \\ 2|j-i}}\gamma(p^{j-i})^2
    \prod_{\substack{i,j \\ 2\nmid j-i}}\gamma(qp^{j-i})}
    {\prod_{i,j} \tau(p^{j-i})} \\
\det T_3 &=
\frac{\prod_{i<j\le 2n}\gamma(p^{j-i})\prod_{\substack{i,j \le 2n+1 \\ 2|j}}
\gamma(qp^{j-i})}
{\prod_{i,j} \gamma(p^{n+j-i})\gamma(p^{n+j+i})}\\
\det T_4 &= \frac{2^n\prod_{i<j \le n} \gamma(p^{2(j-i)})^2
\prod_{\substack{i,j \le 2n+1 \\ 2\nmid i,2|j}} \gamma(qp^{j-i})}
{\prod_{i,j} \tau(p^{j-i})\tau(p^{j+i})}
\end{align*}
\label{th:pq}
\end{theorem}
\begin{proof} Factor exhaustion.  We first view each determinant as a
fractional Laurent polynomial in $q$.  By choosing special values of $q$, we
will find enough factors in each determinant to account for their entire width,
thus determining them up to a rational factor $R(p)$.  (Each determinant is a
centered Laurent polynomial in $q$ with fractional exponents.  The notion of
width make sense for these.)  We will derive this factor by a separate method.

For example, if $0 \le k < n$, then $\det T_1$ is divisible by
$\gamma(qp^{-k})^{n-k}$ because
$$
T_1(p,p^k)_{i,j} = \sum_{\frac{1-k}2\le \ell \le \frac{k-1}2} p^{\ell(n+j-i)}.
$$
Evidently $T_1(p,p^k)$ is a sum of $k$ rank 1 matrices at this specialization,
so its determinant has an $(n-k)$-fold root at $q=p^k$.  Likewise $T(p,p^{-k})$ also has
rank $k$ and $\gamma(qp^k)^{n-k}$ also divides $T_1$.  All four of the
determinants have this behavior.  In each case, the singular values of
$q$ can be read from the product formulas for the determinants.  The only
detail that changes is the form of each rank 1 term, which is summarized in
\tab{t:terms}.

\begin{fulltable}{t:terms}
    {Details of factor exhaustion for \thm{th:pq}}
\begin{tabular}{ccc}
Matrix & Rank 1 terms & Extra $q$ value\\ \hline
$T_1$ & $p^{-\ell i}p^{\ell (n+j)}$ & $\infty$\\
$T_2$ & $p^{-\ell i}p^{\ell j}$ & 1\\
$T_3$ & $(p^{-\ell i} - p^{\ell i})(p^{\ell (n+j)} - p^{-\ell (n+j)})$
& $\infty$\\
$T_4$ & $(p^{-\ell i} - p^{\ell i})(p^{\ell j} - p^{-\ell j})$ &
$\infty$ \\
\end{tabular}
\end{fulltable}

Finally the $q$-independent factor $R(p)$ can be found by examining
the coefficient of the leading power of $q$, or equivalently, taking
the limit $q \to \infty$.  For example
\begin{align*}
\frac{T_1(p,q)_{i,j}}{q^{(n+j-i)/2}} \to \frac1{\gamma(p^{n+j-i})}
    = p^{(n+i+j)/2}C(\vx,\vy)
\end{align*}
as $q \to \infty$ with
$$x_i = p^{-i} \qquad y_j = p^{n+j}.$$
In this case $R(p)$ is given by $\det C_1$ in
\thm{th:cauchy}.  This happens
in each case, although for the matrix $T_2$ it is slightly more convenient
to specialize to $q=1$.  The best extra value of $q$ in all four cases
is given in \tab{t:terms}.
\end{proof}

Finally we evaluate two three-variable Pfaffians which are like
the determinants in \thm{th:pq}.

\begin{theorem}  For $i,j \le 2n$, let
\begin{align*}
T_5(p,q,r)_{i,j} =&\ \frac{\gamma(q^{j-i})\gamma(r^{j-i})}{\gamma(p^{j-i})} \\
T_6(p,q,r)_{i,j} =&\ \gamma(p^{j+i})\gamma(p^{j-i})\biggl(\frac{\gamma(q^{j+i})}{\gamma(p^{j+i})} -
\frac{\gamma(q^{j-i})}{\gamma(p^{j-i})}\biggr) \\
&\  \biggl(\frac{\gamma(r^{j+i})}{\gamma(p^{j+i})} - \frac{\gamma(r^{j-i})}{\gamma(p^{j-i})}\biggr).
\end{align*}
when $i \ne j$ and
\begin{align*}
T_5(p,q,r)_{i,i} &= 0 \\
T_6(p,q,r)_{i,i} &= 0.
\end{align*}
Then
\begin{align*}
\Pf T_5 &= \frac{\prod_{i<j} \gamma(p^{j-i})^4\prod_{i,j}
    \gamma(qp^{j-i})\gamma(rp^{j-i})}{\prod_{i<j \le 2n} \gamma(p^{j-i})} \\
\Pf T_6 &= \frac{\prod_{i<j \le 2n}\gamma(p^{j-i})\prod_{\substack{i,j \le 2n+1 \\ 2|j}}
\gamma(qp^{j-i})\gamma(rp^{j-i})}{\prod_{i<j \le 2n}\gamma(p^{j+i})}.
\end{align*}
\label{th:pqr}
\end{theorem}
\begin{proof} Factor exhaustion in both $q$ and $r$.  If $0 \le k < n$,
then
$$T_5(p,p^k,r)_{i,j} = \sum_{\frac{1-k}2\le \ell \le \frac{k-1}2}
r^{1/2}p^{\ell(j-i)} - \sum_{\frac{1-k}2\le \ell \le \frac{k-1}2}
r^{-1/2}p^{\ell(j-i)}$$
is, as written, a sum of $2k$ rank 1 matrices. Therefore the Pfaffian, whose
square is the determinant, it is divisible by $\gamma(qp^{-k})^{n-k}$.  The
same argument applies to $T_5(p,p^{-k},r)$.  It also applies to $T_5(p,q,p^{\pm k})$ since
$T_5$ is symmetric in $q$ and $r$.

This determines $\Pf T_5$ up to a factor $R(p)$ depending only on $p$.
This factor can be determined by taking the limit $r \to \infty$:
$$\lim_{r \to \infty} \frac{T_5(p,q,r)_{i,j}}{r^{(|i-n-\frac12|+|j-n-\frac12|)/2}} =
\begin{cases} T_1(p,q)_{i,j-n} & i \le n < j \\
-T_1(p,q)_{j,i-n} & j \le n < i \\
0 \end{cases}.$$
In other words, after rescaling rows and columns, $T_5(p,q,r)$
has a block matrix limit:
$$\lim_{r \to \infty} r^{\bullet}T_5(p,q,r) =
\left(\begin{array}{c|c} 0 & T_1(p,q) \\ \hline -T_1(p,q)^T & 0
\end{array}\right).$$
(The bullet $\bullet$ is the exponent above that is different
for different rows and columns.)
This establishes that the leading coefficient of $\Pf T_5(p,q,r)$
as a polynomial in $r$ is
$$(-1)^{\binom{n}{2}} \det T_1(p,q),$$
which in turn determines $R(p)$.

The Pfaffian $\Pf T_6$ is argued the same way.  To find the factor $R(p)$ which
is independent of $q$ and $r$, we take the  limit $r,q \to \infty$.  In this
limit $\Pf T_6$ reduces to a special case of $\Pf C_4$ in \thm{th:cauchy}.
\end{proof}

\begin{remark} Several other specializations of the determinants in
\thm{th:pq} and the Pfaffian in \thm{th:pqr} are
special cases of \thm{th:cauchy} and other determinants and Pfaffians
such as these \cite{Krattenthaler:advanced}:
$$
\det\left\{ x_i^{j-1} \right\} \qquad \det\left\{\gamma(x_i^{2j-1})\right\}
$$
For example, $T_1(q^2,q)$ is also a Cauchy double alternant, while $T_1(p,p^n)$
is the product of two (rescaled) Vandermonde matrices.  Any of these
intersections may be used to determine the $q$-independent factor in the
factor exhaustion method. There are also other round determinants like the ones
in \thm{th:pq} which we do not need, for example
$$
\det\left\{ \frac{\gamma(q^{n+j-i})}{\gamma(p^{n+j-i})} -
\frac{\gamma(q^{n+j+i-1})}{\gamma(p^{n+j+i-1})} \right\}.
$$
These examples suggest the following more general problem: Let $M$ be an $n
\times n$ matrix such that $M_{i,j}$ is a rational polynomial in a fixed
number of variables, such as $p$, $q$, and $r$, and in exponentials of them
such as $p^i$, $q^j$, and $r^n$.  When is $\det M$ round?  What if $M_{ij}$
is a rational polynomial in variables such as $x_i$ and $y_j$?
\end{remark}

\section{Enumerations and divisibilities}
\label{s:prod}

In this section we relate the quantities appearing in the
other sections to prove the results in Section~\ref{s:intro}.

Let
\begin{align}
\vone &= (1,1,1,1,\ldots,1) \nonumber \\
x &= a^2+2+\ba^2 \label{e:xa} \\
y &= \sigma(ba)/\sigma(b\ba) \nonumber \\
z &= \sigma(ca)/\sigma(c\ba). \nonumber
\end{align}
Then most of the generating functions in Section~\ref{s:intro} can
be expressed in terms of the partition functions in Section~\ref{s:det},
\begin{align}
A(n;x) &= \frac{Z(n;\vone,\vone)}{\sigma(a)^{n^2-n}\sigma(a^2)^n}  \label{e:az} \\
A_\HT(2n;x,\pm 1) &= \frac{Z^\pm_\HT(n;\vone,\vone)}
    {\sigma(a)^{2n^2-n}\sigma(a^2)^n} \nonumber \\
A_\U(2n;x,y) &= \frac{Z_\U(n;\vone,\vone)}
    {\sigma(a)^{2n^2-n}\sigma(a^2)^n\sigma(b\ba)^n} \nonumber \\
A_\UU(4n;x,y,z) &= \frac{Z_\UU(n;\vone,\vone)}
    {\sigma(a)^{4n^2-n}\sigma(a^2)^n\sigma(b\ba)^n\sigma(c\ba)^n} \nonumber \\
A_\QT(4n;x) &= \frac{Z_\QT(n;\vone)}
    {\sigma(a)^{4n^2-n}\sigma(a^2)^n} \nonumber \\
A_\O(2n;x) &= \frac{Z_\O(n;\vone)}
    {\sigma(a)^{2n^2-2n}\sigma(a^2)^n} \nonumber \\
A_\UO(8n;x,z) &= \frac{Z_\UO(n;\vone,\vone)}
    {\sigma(a)^{8n^2-3n}\sigma(a^2)^n\sigma(c\ba)^nb^{2n}}, \nonumber
\end{align}
by the definition of the partition functions and the correspondence
between square ice and alternating-sign matrices.
The generating function $A_\OO(4n;x,y)$ requires a slightly different
change of parameters: if $y = b^2/c^2$, then
$$A_\OO(4n;x,y) = \frac{Z_\OO(n;\vone)}
    {\sigma(a)^{4n^2-3n}\sigma(a^2)^nc^{2n}} .$$
The generating function $A_\U(n;x,y)$ is a polynomial of degree $n$ in $y$ and
it is easy to show that the leading and trailing coefficients count VSASMs, so
we can say that
$$A_\V(2n+1;x) = A_\U(n;x,0) = A_\U(n;x,\infty),$$
where by abuse of notation, if $P(x)$ is a polynomial (or a
rational function), $P(\infty)$ denotes the top-degree
coefficient.
Likewise $A_\UU(n;x,y,z)$ has bidegree $(n,n)$ in 
$y$ and $z$ and the corner coefficients count VHSASMs and VHPASMs:
\begin{align*}
A_\VHP(4n+2;x) &= A_\UU(n;x,0,\infty) = A_\UU(n;x,\infty,0)\\
A_\VH(4n+1;x) &= A_\UU(n;x,\infty,\infty) \\
A_\VH(4n+3;x) &= A_\UU(n;x,0,0).
\end{align*}

We can reverse these relations by defining
\begin{align*}
Z^{\pm,(2)}_\HT(n;\vx,\vy) =&\  \frac{\prod_{i,j}\alpha(x_i\by_j)
    (\det M^\pm_\HT)}{\prod_{i<j}\sigma(\bx_ix_j)\sigma(y_i\by_j)}
     \\
Z^{(2)}_\UU(n;\vx,\vy) =&\ \frac{\prod_{i,j}\alpha(x_i\by_j)\alpha(x_iy_j)
    (\det M_\UU)}{\prod_{i<j}\sigma(\bx_ix_j)\sigma(y_i\by_j)
    \prod_{i\le j}\sigma(\bx_i\bx_j)\sigma(y_iy_j)} \\
Z^{(k)}_\QT(n;\vx) =&\  \frac{\prod_{i<j \le 2n}\alpha(\bx_ix_j)
    (\Pf M^{(k)}_\QT)}{\prod_{i<j \le 2n} \sigma(\bx_ix_j)} \\
Z^{(2)}_\OO(n;\vx) =&\ \frac{\prod_{i<j \le 2n}\alpha(x_ix_j)
    (\Pf M_\OO)}{\prod_{i<j \le 2n} \sigma(\bx_ix_j)} \\
Z^{(k)}_\UO(n;\vx) =&\ \frac{\prod_{i<j\le 2n}\alpha(\bx_i x_j)\alpha(x_i x_j)(\Pf M^{(k)}_\UO)}
    {\prod_{i < j \le 2n} \sigma(\bx_i x_j)\prod_{i\le j\le 2n}\sigma(x_i x_j)}
\end{align*}
and 
\begin{align*}
A^{(2)}_\HT(2n;x,\pm 1) &= \sigma(a)^{n-n^2} Z^{\pm,(2)}_\HT(2n;\vone,\vone) \\
A^{(2)}_\UU(4n;x,y,z)   &= \sigma(a)^{n-2n^2}\sigma(b\ba)^{-n}\sigma(c\ba)^{-n}
    Z^{(2)}_\UU(4n;\vone,\vone) \\
A^{(k)}_\QT(4n;x)       &= \sigma(a)^{n-2n^2} Z^{(k)}_\QT(4n;\vone) \\
A^{(2)}_\OO(4n;x,y)     &= c^{-2n} \sigma(a)^{2n-2n^2} Z^{(2)}_\OO(4n;\vone) \\
A^{(1)}_\UO(8n;x,z)     &= \sigma(a)^{n-4n^2} Z^{(1)}_\UO(8n;\vone,\vone) \\
A^{(2)}_\UO(8n;x,z)     &= \sigma(a)^{n-4n^2}\sigma(c\ba)^{-n} Z^{(2)}_\UO(8n;\vone,\vone).
\end{align*}
using the same correspondence between $a$, $b$, and $c$ with $x$, $y$,
and $z$ (which slightly different in the case of OOSASMs).
The observation that all of these quantities must be polynomials
establishes the factorizations of $A_\HT$, $A_\OO$, $A_\QT$,
$A_\UU$, and $A_\UO$ in \thm{th:factor}.

For vectors
$\vx$ and $\vy$, let
$$(\vx,\vy) = (x_1,\ldots,x_n,y_1,\ldots,y_n)$$
denote their concatenation, and let exponentiation of vectors
denote coordinate-wise exponentiation:
$$\vx^k = (x_1^k,x_2^k,\ldots,x^k_n).$$
Then the matrices
\begin{gather*}
M(2n;(\vx,\vx^{-1}),(\vy,\vy^{-1})) \\
M_\HT(2n;(\vx,\vx^{-1}),(\vy,\vy^{-1}))
\end{gather*}
commutes with the permutation matrix
$$\left(\begin{array}{c|c} 0 & I_n \\ \hline I_n & 0 \end{array}\right),$$
where $I_n$ is the $n \times n$ identity matrix.  Similarly
the matrices
\begin{gather*}
M(2n+1;(\vx,1,\vx^{-1}),(\vy,1,\vy^{-1})) \\
M_\HT(2n+1;(\vx,1,\vx^{-1}),(\vy,1,\vy^{-1}))
\end{gather*}
commute with
$$P = \left(\begin{array}{c|c|c} 0 & 0 & I_n \\
\hline 0 & 1 & 0 \\ \hline I_n & 0 & 0 \end{array}\right).$$
If $\vx$ has length $2n$ and $b^2 = -c^2$, then
$$M_\OO(4n;(\vx,\vx^{-1}))$$
commutes with
$$P = \left(\begin{array}{c|c} 0 & I_{2n} \\ \hline I_{2n} & 0 \end{array}\right).$$
In each case we can restrict decompose $M$, $M_\HT$, and $M_\OO$ into
blocks corresponding to the eigenspaces of $P$.  The block with
eigenvalue $-1$ is in the three cases equal to $M_\U$ and proportional
to $M_\UU$ (with $b=c=\ii$) and $M^{(1)}_\UO$.  This results in
the factorizations of $A(n;x)$, $A^{(2)}_\HT(n;x,1)$, and 
$A^{(2)}_\OO(4n;x,-1)$ in \thm{th:factor}.  Another
factorization in \thm{th:factor} is that of $A^{(2)}_\HT(2n;x,-1)$.
To establish this, observe that the matrix $M^-_\HT(2n;\vx,\vx)$
is antisymmetric, and that it is proportional to $M^{(1)}_\QT(2n;\vx)$.
The latter matrix is employed for its Pfaffian, while the former
for its determinant, which is the square of the Pfaffian.

The final case of \thm{th:factor} is the relation between $A_\V(2n+1;x)$
and $A_\U(2n;x)$.  This relation is established by observing 
that $b$ appears only in the normalization factor for $Z_\U(n;x)$
and not in the matrix $M_\U(n;x)$; the only step is to change variables
from $b$ to $y$.

\begin{fulltable}{t:spec}{Specializations of partition function determinants
    and Pfaffians}
\begin{tabular}{l@{\hspace{.4cm}}l@{\hspace{.2cm}}l@{\hspace{.2cm}}l}
Enumeration & Parameters   & Section~\ref{s:det}     & Section~\ref{s:factor} \\ \hline
$A(n;1)$                & $a=\omega_3$      & $M(\vq,\vq(n))$       & $T_1(q^3,q)$ \\
$A(n;2)$                & $a=\omega_4$      & $M(\vx,\vy)$          & $C_1(\vx^2,\vy^2)$ \\
$A(n;3)$                & $a=\omega_6$      & $M(\vq,\vq)$          & $T_2(q^3,q)$ \\
$A_\HT(2n;1,1)$         & $a=\omega_3$      & $M^+_\HT(\vq,\vq(n))$ & $T_1(q^3,q^2)$ \\
$A_\HT(2n;2,1)$         & $a=\omega_4$      & $M^+_\HT(\vq,\vq)$    & $T_2(q^2,q)$ \\
$A_\V(2n+1;1)$          & $a=\omega_3$      & $M_\U(\vq,\vq(n))$    & $T_3(q^3,q)$ \\
$A_\V(2n+1;2)$          & $a=\omega_4$      & $M_\U(\vx,\vy)$       & $C_2(\vx^2,\vy^2)$ \\
$A_\V(2n+1;3)$          & $a=\omega_6$      & $M_\U(\vq,\vq)$       & $T_4(q^3,q)$ \\
$A^{(2)}_\UU(4n;1,1,1)$ & $a=\omega_3,b=c=\omega_4$ & $M_\UU(\vq,\vq(n))$ & $T_3(q^3,q^2)$ \\
$A^{(2)}_\UU(4n;2,1,1)$ & $a=b=c=\omega_4$ & $M_\UU(\vq,\vq)$      & $T_4(q^2,q)$ \\
$A^{(2)}_\VHP(4n;1)$    & $a=b=\omega_3,c=\ba$ & $M_\UU(\vq,\vq(n))$  & $T_3(q^3,q)$ \\
$A^{(1)}_\QT(4n;1)$     & $a=\omega_3$      & $M^{(1)}_\QT(\vq)$    & $T_5(q^3,q,q)$ \\
$A^{(1)}_\QT(4n;2)$     & $a=\omega_4$      & $M^{(1)}_\QT(\vq)$    & $T_5(q^4,q^2,q)$ \\
$A^{(1)}_\QT(4n;3)$     & $a=\omega_6$      & $M^{(1)}_\QT(\vq)$    & $T_5(q^6,q^3,q^2)$ \\
$A^{(2)}_\QT(4n;1)$     & $a=\omega_3$      & $M^{(2)}_\QT(\vq)$    & $T_5(q^3,q^2,q)$ \\
$A^{(2)}_\QT(4n;2)$     & $a=\omega_4$      & $M^{(2)}_\QT(\vx)$    & $C_3(\vx^2)$ \\
$A_\O(2n;1)$            & $a=\omega_3$      & $M_\O(\vq)$           & $T_6(q^3,q,\infty)$ \\
$A^{(1)}_\UO(8n;1)$     & $a=\omega_3,c=\omega_4$ & $M^{(1)}_\UO(\vq)$ & $T_6(q^3,q,q)$ \\
$A^{(2)}_\UO(8n;1,1)$   & $a=\omega_3,c=\omega_4$ & $M^{(2)}_\UO(\vq)$ & $T_6(q^3,q^2,q)$
\end{tabular}\eatline\end{fulltable}

Finally we establish the round enumerations in Theorems~\ref{th:main},
\ref{th:3main}, and \ref{th:extra}.  We review the argument from
Reference~\citealp{Kuperberg:asm} for $A(n)$, which is an illustrative case.
Let $\omega_n = \exp(\pi\ii/n)$, where $\ii^2 = -1$.  Equation~\ref{e:xa}
implies the following correspondence between $x$ and $a$:
$$\begin{array}{c@{\quad}c@{\quad}c}
a = \omega_3 & \implies & x = 1 \\
a = \omega_4 & \implies & x = 2 \\
a = \omega_6 & \implies & x = 3.
\end{array}$$
For any of these values of $x$ or $a$, we would like to evaluate
the partition function $Z(n;\vone,\vone)$ to find $A(n;x)$
by equation~\ref{e:az}.  Unfortunately the matrix $M(n;\vone,\vone)$
is singular.  So instead we will find its determinant along
a curve of parameters that includes $(\vone,\vone)$.  More precisely,
let 
$$\vq(k) = (q^{(k+1)/2},q^{(k+2)/2},\ldots,q^{(k+n)/2})$$
and
$$\vq = \vq(0).$$
Then 
$$\lim_{q \to 1} \vq(k) = \vone,$$
and if we assume equation~\ref{e:xa},
$$M(n;\vq(0),\vq(k)) = \frac1{-q^{k+j-i}+x-2-q^{i-j-k}}.$$
If we further set $x=1$ and $k=n$, then 
$$M(n;\vq(0),\vq(n)) = -T_1(q^3,q)$$
has a round determinant by \thm{th:pq}.  Computing $A(n;x)$ is then a
routine but tedious simplification of round products. The argument for most of
the other enumerations is the same, except that the curve of parameters is
$(\vq,\vq(n))$ for 1-enumeration in the determinant cases, $(\vq,\vq)$  for 2-
and 3-enumeration in the determinant cases, and $\vq$ in the Pfaffian cases. 
Each of the matrices is then proportional to some matrix $T_i$ from
Theorem~\ref{th:pq} or \ref{th:pqr}.  The
determinants and Pfaffian for three of the 2-enumerations are round without
specializing the parameters and instead reduce to a matrix $C_i$
from \thm{th:cauchy}.  These variations of the argument are
summarized in \tab{t:spec}.

\section{Discussion}
\label{s:discuss}

\begin{table*}
\begin{center}
\begin{tabular}{l@{\hspace{.5cm}}l@{\hspace{.5cm}}l@{\hspace{.5cm}}l}
Factor                 & $n=1$ & $n=2$             & $n=3$ \\ \hline
$A_\V(2n+1;x)$         & 1     & $x+2$             & $x^3+6x^2+13x+6$ \\
$\tA_\V(2n;x)$         & 1     & $x+6$             & $x^3+12x^2+70x+60$ \\
$A^{(1)}_\QT(4n;x)$    & 1     & $x+3$             & $x^3+8x^2+25x+15$ \\
$A^{(2)}_\QT(4n;x)$    & 2     & $6x+4$            & $20x^3+60x^2+52x+8$ \\
$A^{(2)}_\UU(4n;x,1,1)$    & $x+4$ & $x^3+9x^2+40x+16$ & $x^6+16x^5+125x^4+629x^3+1036x^2+560x+64$ \\
$\tA^{(2)}_\UU(4n;x)$  & 1     & $x+1$             & $x^3+4x^2+5x+1$ \\
$A^{(2)}_\VHP(4n+2;x)$ & 1     & $2x+1$            & $5x^3+12x^2+8x+1$ \\
$A_\O(2n;x)$           & 1     & 3                 & $x^2+10x+15$ \\
$A^{(2)}_\OO(4n;x,1)$  & 2     & $2x+4$            & $2x^3+20x^2+28x+8$ \\
$A^{(1)}_\UO(8n;x)$    & 1     & $x^2+5x+3$        & $x^6+14x^5+82x^4+210x^3+239x^2+115x+15$ \\
$\tA^{(1)}_\UO(8n;x)$  & 2     & $2x^2+20x+20$     & $2x^6+42x^5+420x^4+1680x^3+2892x^2+2040x+360$ \\
$A^{(2)}_\UO(8n;x,1)$    & $x+4$ & $x^4+12x^3+65x^2+104x+16$ & $x^9+24x^8+275x^7+1966x^6+8215x^5+19144x^4$ \\
& & & $+21777x^3+10028x^2+1712x+64$
\end{tabular}
\end{center}
\caption{Irreducible $x$-enumerations}
\label{t:irred}
\end{table*}

Even though each section of this article considers many types of
alternating-sign matrices or determinants in parallel, the  work of enumerating
symmetry classes of ASMs is far from finished.  Robbins \cite{Robbins:symmetry}
conjectures formulas for the number of VHSASMs and for the number of odd-order
HTSASMs, QTSASMs, and DASASMs in addition to the enumerations that we have
proven.  \thm{th:z} yields a determinant formula for the number of VHSASMs,
obtained from the more general partition function $Z_\UU(n;\vx,\vy)$ by setting
$a = \omega_3$ and $b = c = \omega_3^{\pm 1}$, and all other parameters to 1). 
In the enumeration of $4n+1 \times 4n+1$ VHSASMs, where $b=c=\omega_3$,
experiments indicate that
$$\det M_\UU(\vq(-\frac12),\vq(n-\frac12))$$
is round, but we cannot prove this.  In the other case, $4n+3 \times 4n+3$
VHSASMs, where $b=c=\omega^{-1}$, we could not even find a curve for which the
determinant is round.  This strange behavior of VHSASMs is one illustration
that although we have put many classes of ASMs under one roof, the
house is not completely in order.

For the other three classes conjecturally enumerated by Robbins, we could not
even find a determinant formula.  Nonetheless we conjecture:

\begin{question} Can DSASMs, DASASMs, TSASMs, and odd-order HTSASMs
and QTSASMs can be $x$-enumerated in polynomial time?
\end{question}

The polynomials listed in \tab{t:irred} appear to be generating functions of
some type, but in most cases there is not even a proof that their coefficients
are non-negative.  (We have more data than is shown in the table; the
multivariate polynomials $A^{(2)}_\UU(4n;x,y,z)$, $A^{(2)}_\OO(4n;x,y)$, and
$A^{(2)}_\UO(4n;x,y)$ also appear to be non-negative.) Some of them are
conjecturally related to cyclically symmetric plane partitions
\cite{Robbins:symmetry}. Indeed they are related to each other in strange ways.
For example \thm{th:extra} establishes that if we take $x=1$, three of the
polynomial series ($A_\V$, $A_\O$, and $A^{(2)}_\VHP$) become equal, as if to
suggest that VSASMs can be $x$-enumerated in three different ways!

\begin{question} Do the polynomials in \tab{t:irred}
$x$-enumerate classes of alternating-sign matrices?
\end{question}

OSASMs include the set of off-diagonal permutation matrices, which can be
interpreted as the index set for the usual combinatorial formula for the
Pfaffian.  Like ASMs, their number is round.  These observations, together with
the known formulas due to Mills, Robbins, and Rumsey \cite{MRR:asm} motivate
the following question:

\begin{question} Are there formulas for the Pfaffian of a matrix involving
OSASMs that generalize the determinant formulas involving ASMs?
\end{question}

Neither any of the enumerations that we establish, nor the various
equinumerations that they imply, have known bijective proofs. Nor is it even
known that two equinumerous types of ASMs index bases of the same vector
space.  For example, can one find an explicit isomorphism between the vector
space of formal linear combinations of $2n \times 2n$ OSASMs and the vector
space of formal linear combinations of $2n+1 \times 2n+1$ VSASMs?

Sogo found that $Z(n;\vone,\vone)$ satisfies the Toda chain (or Toda molecule)
differential hierarchy~\cite{Sogo:orthogonal,KZ:domain}.

\begin{question} If $\vx$ and $\vy$ are set to $\vone$, do the partition
functions in \thm{th:z} and \tab{t:irred} satisfy natural differential
hierarchies?
\end{question}

Many other solutions to the Yang-Baxter equation are known
\cite{Drinfeld:quantum}.  The six-vertex solution corresponds to the Lie
algebra $\sl(2)$ together with its 2-dimensional representation;
there are solutions for other simple Lie algebras and their representations.

\begin{question} Do square ice and Izergin-Korepin-type determinants generalize
to other solutions of the Yang-Baxter equation? \end{question}

Although our simultaneous treatment of several classes of ASMs is not
especially short, the argument for any one alone is relatively simple.  We
speculate that the methods of Lagrange interpolation (used in
Section~\ref{s:det}) and factor exhaustion (the topic of
Section~\ref{s:factor}) simplify many proofs of product formulas.  I. J. Good's
short proof of Dyson's conjecture \cite{good:dyson} also uses Lagrange
interpolation.


\providecommand{\bysame}{\leavevmode\hbox to3em{\hrulefill}\thinspace}

\end{document}